\documentclass[a4paper, 11pt]{article}

\usepackage[english]{babel}
\usepackage{authblk}

\usepackage{braket}
\usepackage{boxedminipage}
\usepackage{bbold}
\usepackage{titlesec}
\usepackage{amsmath}
\usepackage{fancybox,xcolor}
\usepackage{empheq}
\usepackage{soul}

\definecolor{shadowcolor}{RGB}{0, 0, 102}

\usepackage{shadow}
\usepackage{amsthm}
\usepackage{amssymb}
\usepackage{graphicx}
\usepackage[colorlinks=true, allcolors=blue]{hyperref}

\usepackage{amssymb,empheq,fancybox}

\newsavebox{\mybox}

\usepackage{geometry}
\usepackage{stmaryrd}
\usepackage{amsmath}
\usepackage{graphicx}
\usepackage[colorlinks=true, allcolors=blue]{hyperref}

\newcommand{\Da}{ {\cal D }}
\newcommand{\Pa}{ {\cal P }}
\newcommand{\Va}{ {\cal V }}
\newcommand{\Qa}{ {\cal Q }}
\newcommand{\La}{ {\cal L }}

\newcommand{\Ga}{ {\cal G }}
\newcommand{\Ma}{ {\cal M }}
\newcommand{\Na}{ {\cal N }}
\newcommand{\Ea}{ {\cal E }}

\newcommand{\Ca}{ {\cal C }}
\newcommand{\Ha}{ {\cal H }}
\newcommand{\Sa}{ {\cal S }}

\newcommand{\E}{ {\mathbb{E}}}
\newcommand{\Prb}{ {\mathbb{P}}}
\newcommand{\R}{ {\mathbb{R}}}
\newcommand{\N}{ {\mathbb{N}}}

\newcommand{\Xa}{ {\cal X }}

\newcommand{\PP}{\mathbb{P}}
\newcommand{\NN}{\mathbb{N}}
\newcommand{\EE}{\mathbb{E}}
\newcommand{\point}{\mbox{\LARGE .}}
\newcommand{\CC}{\mathbb{C}}

\newcommand{\Pbar}{\chi}

\newtheorem{thm}{Theorem}
\newtheorem{lmm}{Lemma}
\newtheorem{cor}{Corollary}

\newtheorem{prop}{Proposition}

\begin{document}

\title{On the Mathematical foundations of Diffusion Monte~Carlo}
\author{Michel Caffarel$^{1}$, Pierre Del Moral$^{2}$ and Luc de Montella$^{2,3}$}
\affil{\footnotesize $^{1}$University of Toulouse, Lab. de Chimie et Physique Quantiques, Toulouse , 31062, FR. {\footnotesize E-Mail:\,} \texttt{\footnotesize caffarel@irsamc.ups-tlse.fr}}
\affil{\footnotesize $^{2}$Centre de Recherche Inria Bordeaux Sud-Ouest, Talence, 33405, FR. {\footnotesize E-Mail:\,} \texttt{\footnotesize pierre.del-moral@inria.fr,luc.de-montella@inria.fr}}
\affil{\footnotesize $^{3}$Naval Group,  Bouguenais, 44340, FR. {\footnotesize E-Mail:\,} \texttt{\footnotesize luc.demontella@naval-group.com}}
\maketitle
\begin{abstract}
	The Diffusion  Monte Carlo method with constant number of walkers, also called Stochastic Reconfiguration as well as Sequential Monte Carlo, is a widely used  Monte Carlo methodology for computing the ground-state energy and wave function of quantum systems.  In this study, we present the first mathematically rigorous analysis of this class of stochastic methods on non necessarily compact state spaces, including linear diffusions evolving in quadratic absorbing
potentials, yielding what seems to be the first result of this type for this class of models. 
We present a novel and general mathematical framework with easily checked Lyapunov stability conditions  that ensure the uniform-in-time convergence of Diffusion Monte Carlo estimates towards the top of the spectrum of Schrödinger operators. For transient free evolutions, we also present a divergence blow up of the estimates w.r.t. the time horizon even when the asymptotic fluctuation variances are uniformly bounded.
We also illustrate the impact of these results in the context of generalized coupled quantum harmonic oscillators with non necessarily reversible nor stable diffusive particle and a quadratic energy absorbing well associated with a semi-definite positive matrix force.\end{abstract}


\section{Introduction}
The many-body Schrödinger equation describes interacting quantum particles. Depending on the domain of application,
these particles may represent electrons in solid-state physics or quantum chemistry, 
nucleons in nuclear physics, atoms in quantum liquid physics,
or coupled modes of oscillators in molecular spectroscopy, among the main applications.
Except for trivial quantum systems, it is impossible to solve this equation analytically.
The diffusion  Monte Carlo method (abbreviated DMC) provides a powerful stochastic
approach to numerically approximate the ground state energy and wave function of Schrödinger operators.

The DMC methodology has a long and rich history, dating back to its first mention in 1949 
by Ulam and Metropolis in \cite{metropolis1949monte}. The idea was first implemented 
by Donsker and Kac \cite{donsker1950sampling}, and by Kalos~\cite{kalos1962monte} in the early 1960s.  
Over the years, the physics community has proposed numerous variants of Diffusion Monte Carlo, known by various names such as
Green's function Monte Carlo,\cite{kalos1962monte,binder} Fixed-Node Diffusion Monte Carlo,\cite{reynolds1982}, Pure Diffusion Monte Carlo,
\cite{caffarel1986treatment,caffarel1988development}
Stochastic Reconfiguration Monte
Carlo,\cite{hetherington1984observations,sorella1998green,sorella2000green,assaraf2000diffusion} 
and Reptation Monte Carlo,\cite{moroni1999} to cite the main ones.
Despite their apparent diversity, all these approaches are fundamentally based, in one way or another, 
on the stochastic simulation of a specific implementation of the Feynman-Kac formula with importance sampling.

For a more detailed discussion on the origins and the applications of these Monte Carlo techniques 
in physics we refer the reader to the recent review article \cite{mareschal2021early} 
as well as to \cite{ma2023ground,foulkes2001quantum} and references therein.

The version of interest employed here is the DMC method with a fixed number of walkers, 
commonly known in physics as \textit{Stochastic Reconfiguration Monte Carlo}; see the pioneering article 
by Hetherington~\cite{hetherington1984observations}, followed by Sorella and 
co-authors~\cite{sorella1998green,sorella2000green} and by the first author 
and his co-workers in \cite{assaraf2000diffusion}.

In mathematics, the methodology may also be referred to by different names, 
such as genetic algorithm with selection and mutation, population Monte Carlo 
or sequential Monte Carlo~\cite{del2003particle,del2004particle,Rousset2006OnTC,articleCancEric,el2007diffusion}. 
For a more thorough discussion on these application model areas we refer to the 
books \cite{bookFK,del2013mean} and references therein.

These sequential Monte Carlo methods do not rely on biased variational techniques. 
They can be seen as a sophisticated genetic-type Monte Carlo methodology to simulate interacting quantum 
many-body systems. Various asymptotic results have been derived, including central limit theorems 
and large deviation principles, see for instance  \cite{del1998large,dawson2005large} and 
\cite{del1999central,chopin2004central}, as well as the books \cite{DelMoral2000,bookFK,del2013mean} for an overview.

Our work concerns less studied non-asymptotic and time-uniform problems. 
Recalling that the estimation of ground state energies relies on the limiting behavior of the walkers' evolution 
in the DMC method, it is therefore crucial to obtain uniform-in-time convergence estimates. 
Despite its importance, there is a notable gap in the literature and very few results have been proven 
in this respect. To the best of our knowledge, such uniform controls are mainly valid for compact state 
space models, see for instance \cite{del2001stability,bookFK,del2013mean} 
as well as \cite{whiteley2012sequential}. Surprisingly, the theoretical efficiency of the DMC method 
has never been verified rigorously even in basic linear-Gaussian scenarios such as the simple and well known harmonic oscillator.
In this paper, we address this gap by establishing the first uniform-in-time convergence estimates 
that apply to general state space models including the coupled harmonic oscillators presented in~\cite{del2023coupled}.

Our approach is partly based on recent developments on the stability of positive semigroups presented in \cite{arnaudon2023lyapunov,del2021stability}, see also the analysis of generalized coupled harmonic oscillators presented in~\cite{del2023coupled}.
In the present article, we provide a natural Lyapunov condition that ensures the exponential stability  of possibly time varying positive semigroups on non necessarily compact state spaces (cf. (\ref{GQV}) and the local conditions (\ref{min-cond-Q})). In the context of time homogeneous positive semigroups, these conditions ensure the existence of an unique leading eigen-triple (see for instance (\ref{triplet-cond})). We underline that these results  do not rely on any reversibility-type condition, nor on some spectral theorem. They can be seen as an extended version of Perron-Frobenius and Krein-Rutman theorems for possibly time varying positive operators.

We present a nonlinear Markov chain interpretation of the DMC methodologies. In this interpretation, the genetic type evolution of the walkers  can be seen as a mean field particle simulation of a nonlinear Markov chain (see Section~\ref{subsection:DMC}).
In this context, we present an auxiliary Lyapunov condition that depends on the potential function and the free evolution of the walkers that ensures the time uniform performance of the DMC methodology (cf. condition (\ref{GQV}),  as well as Theorem~\ref{thm:optimal_predictor_unif_bound} and Corollary~\ref{cor:control_limit}). 
We illustrate this condition in the context of generalized coupled harmonic oscillators for a linear diffusive-type particle and a quadratic energy absorbing well associated with a semi-definite positive matrix force.  In this context, we also show that the DMC methodology may  diverge  when the free evolution of the walkers is unstable {\em  for any fixed number of walkers, even if the asymptotic variance of the Central Limit Theorem is uniformly bounded} with respect to the time parameter (see Proposition~\ref{thm:div_oscillo}, as well as Section~\ref{subsection:hamo_1_d} and Proposition~\ref{prop-fi}).

In the context of absorbing wells centered at the origin, this study leads us to conjecture that stable free evolution transitions is a  necessary and sufficient condition for the DMC method to be 
uniformly convergent w.r.t. the time horizon.

Additionally, we propose and, to some extent, establish the validity of an importance sampling transformation 
to overcome this difficulty. This type of technique is
related but not identical to the use of guiding wave 
functions in physics to direct the Monte Carlo moves to improve the efficiency 
of the DMC method \cite{PhysRev.128.1791, PhysRevA.9.2178}. In contrast with conventional guiding waves techniques our approach is based on conditional free evolutions transitions and survival weight potential functions (cf. Section~\ref{subsection:DMC} and Section~\ref{subsection:importance sampling}).

The rest of the article is organized as follows: In Section \ref{section:description_model}, we provide a detailed 
description of the general framework in which our study is set, as well as the theoretical 
foundations on which our proof will be based.

Section \ref{section:main_result} is devoted to the presentation of our main results. 
Section \ref{section:proof_bound} is mainly concerned the detailed proofs of time-uniform estimates.

Section \ref{section:harmonic_oscillator} is devoted to the application of our convergence result to generalized 
coupled harmonic oscillators \cite{del2023coupled,messiah1961quantum, griffiths_introduction_2018}. These models arise 
in various fields such as molecular spectroscopy\cite{herzberg1945}, quantum optics \cite{PhysRevE.102.052213}, quantum cryptography \cite{PhysRevLett.67.661} 
and photosynthesis \cite{articlePhotosynthesis}.  In signal processing, the harmonic oscillator and the DMC methods coincides with the Kalman and the particle filter \cite{articleIPS,Klmn1961NewRI}.

\section{Description of the models}\label{section:description_model}

\subsection{Free evolution semigroups}\label{subsection:markov}
Consider a Markov chain $X_n$ indexed by $n\in\NN$ and taking values in a locally compact Polish space
$(E,\Ea)$, where $\Ea$ is the Borel $\sigma$-field on $E$. 
Let $\Ca(E)$ be the algebra of continuous measurable functions on $E$. We also define
 $\Ca_b(E) \subset \Ca(E)$ as the sub-algebra of bounded measurable continuous functions endowed with the supremum  norm $\|\point \|$.  
With a slight abuse of notation, we denote by $0$ and $1$ the null and unit scalars as well as the null and unit functions on $E$ and we denote by $I:x\in E\mapsto I(x)=x$ the identity function on $E$.  Throughout, we will use \(a\), \(b\), \(c\), \(c_1\), or \(c_2\) to represent positive constants, whose values may vary from line to line.

For $n\in\N^*$, we consider the Markov transitions $P_{n}$ associated with $X_n$, and assume that they are Feller; in the sense that for any $f\in\Ca_b(E)$ we have $P_{n}(f)\in \Ca_b(E)$, with the function $P_{n}(f)$ defined for any $x\in E$ by the integral operator

\begin{equation*}
P_{n}(f)(x) ~:=~\int_E P_{n}(x, dy)f(y) ~=~ \E(f(X_n) ~|~ X_{n-1}=x). 
\end{equation*}

Let $\Ca_{\infty}(E)\subset \Ca(E)$ be the sub-algebra of  uniformly positive continuous functions $V$ that grow at infinity; that is, 
 for any $r\geq V_{\star}:=\inf_E V>0$,  the $r$-sub-level set $\Va(r):=\{V\leq r\}\subset E$ is a non-empty compact subset. 
 We further assume that there exists a  $P$-Lyapunov function $V\in \Ca_{\infty}(E)$; in the sense that $V(E)\subset[1, \infty)$ and there exists $\epsilon\in [0,1)$ and $c\in\R$ such that for any $n\in\N^*$ we have

\begin{equation} \label{lyap-cond}
 P_{n}(V)\leq \epsilon V+c.
\end{equation}

Let $\Ca_V(E)\subset \Ca(E)$ be the sub-space of functions $f\in \Ca(E)$ such that $f/V$ is bounded, equipped with the norm
$\Vert f\Vert_V:=\Vert f/V\Vert$. The Markov semigroup associated with the Markov chain $X_n$  is defined for any $f\in\Ca_V(E)$ by

\begin{equation*}
P_{k, n}(f)(x)  ~:=~ \E(f(X_n) ~|~ X_{k}=x).
\end{equation*}
Condition \eqref{lyap-cond} ensures that  $P_{k, n}$ is
 $V$-Feller  in the sense that for $f\in \Ca_V(E) $ we have
$P_{k,  n}(f)\in \Ca_V(E)$. 
To ensure the semigroup $P_{k,  n}$ is exponentially stable~\cite{arnaudon2023lyapunov}, we assume  the integral operator $$P_n(x,dy)=p_n(x,y)\nu(dy)$$ has a density $p_n$ w.r.t. some Radon measure
$\nu$ satisfying for some $r_1>0$ and for any $r\geq r_1$ the local minorization condition 
\begin{equation}\label{min-cond}
\left\{
    \begin{array}{ll}
        0<\inf_{n\in\N^*}\inf_{\Va(r)^2} p_n\leq \sup_{n\in\N^*}\sup_{\Va(r)^2} p_n < \infty \\ \\
        0<\nu(\Va(r))<\infty.
    \end{array}
\right.
\end{equation}

The $V$-norm semigroup contraction techniques developed in Section 8.2 in~\cite{delmoral:halartilce1}
(see also Lemma 2.3 in~\cite{del2021stability} and Theorem 2.2 in \cite{arnaudon2023lyapunov}), ensure that for any $\mu\in \Pa_V(E)$, there exists some parameters $a,b>0$ such that for any $k\leq n$, and any $\mu_1,\mu_2\in\Pa_V(E)$ we have

\begin{equation}
\left |\left |\mu_1 P_{k,n} - \mu_2 P_{k, n} \right |\right |_{V} \leq ae^{-b(n-k)}||\mu_1 - \mu_2||_{V} \label{contraction-P}.
\end{equation}

Note that the r.h.s condition in \eqref{min-cond} is met as soon as $V$ has compact sub-level sets with non empty interior and $\nu$ is a Radon measure with full support; that is $\nu$ is finite on compact sets and strictly
positive on non-empty open sets. For time-homogeneous models, the l.h.s. minorization condition is satisfied as soon as $(x,y)\in (E^{\circ})^2\mapsto p_n(x,y)$ is a continuous positive function on the interior $E^{\circ}$ of the set  $E$.

\subsection{Feynman-Kac semigroups}\label{subsection:markov_transition}

 We denote by $\Ca_0(E):=\{1/V~:~V\in \Ca_{\infty}(E)\}\subset \Ca_b(E)$ the sub-algebra of bounded continuous positive functions $h$ that vanish at infinity;  that is, for any $0<\epsilon\leq \Vert h\Vert<\infty$ the $\epsilon$-super-level set $\{h\geq \epsilon\}\subset E$ is a non empty compact subset.

Consider a family of strictly positive bounded continuous functions $G_n$  such that

\begin{equation}
{G_n P_{n+1}(V)}/{V} \in \Ca_0(E) \label{cond_G_P(V)/V}
\end{equation}

Condition \eqref{lyap-cond} ensures that this is met for any family of functions in $\Ca_0(E)$.

We associate with it the discrete time generation 
Feynman-Kac semigroups

\begin{equation*}
        Q_{k, n}(f)(x) = \E\left (f(X_n)\prod\limits_{p=k}^{n-1}G_p(X_p) ~|~X_k = x\right ) 
\end{equation*}
and
\begin{equation*}
       \widehat Q_{k, n}(f) :=  G_{k}^{-1} Q_{k, n}(G_nf),
\end{equation*}
with $ G_{k}^{-1}:=1/ G_{k}$. Note that condition \eqref{cond_G_P(V)/V} implies that for any positive function $f\in \Ca_V(E)$ we have $Q_{n, n+1}(f)/V\leq Q_{n, k}(V)/V \in \Ca_0(E)$. Therefore, $Q_{n, n+1}$ is $V$-Feller. In the same vain, we check that $Q_{k,n}$ is $V$-Feller.

To simplify notation, for $k=(n-1)$ we sometimes write $Q_n$ and $\widehat{Q}_n$  instead of $Q_{n-1, n}$ and $\widehat{Q}_{n-1, n}$. In this notation, we have
$$
Q_n(f)(x)=G_{n-1}(x)~P_n(f)(x)\quad \mbox{\rm and}\quad
\widehat{Q}_n(f)(x)=P_n(G_nf)(x).
$$
We also use the convention $Q_{n,n}=P_{n,n}=Id$, the identity operator.

Let  $\Ma_b(E)$ be the set of bounded signed measures on $E$. Also, let $\Pa(E) \subset \Ma_b(E)$ be the convex subset of probability measures on $E$ and denote by $\Pa_V (E)$ the convex set of probability measures $\mu \in \Pa(E)$ such that $\mu(V) < +\infty$.
The left action of $Q_n$ on $\Pa_V(E)$ is given for any $(\eta,f) \in (\Pa_V(E),\Ca_V(E))$ by the formula 
\begin{align}
(\eta Q_n)(f) :=\eta(Q_n(f)) &= \int \eta(dx)Q_n(f)(x) \notag \\
    &= \int \eta(dx)Q_n(x, dy)f(y). \label{action_measure}
\end{align}
By Fubini's theorem, the integration order doesn't matter. Thus to simplify notation, we sometimes write $\eta Q_n(f)$ instead of $(\eta Q_n)(f)$ or $\eta (Q_n(f))$.

We further assume the Lyapunov function $V$ introduced in \eqref{lyap-cond} is a $Q$-Lyapunov function in the sense that \eqref{lyap-cond} holds and there exists $\Theta \in \Ca_0(E)$ and a compact subset  $K\subset E$ such that for any $n\geq 1$  we have
\begin{equation}\label{GQV}
\left\{
    \begin{array}{ll}
        Q_{n}(V)/V \leq \Theta \\ \\
        (G_{n-1}(x) - G_{n-1}(y))(\mathbb{1}_{E\setminus K}(x)V(x) - \mathbb{1}_{E\setminus K}(y)V(y)) \leq 0.
    \end{array}
\right.
\end{equation}
Note that the second condition in \eqref{GQV} holds as soon as there exists $G\in\Ca_0(E)$ such that for any $G_n \leq G$, for any $n\geq 0$. This condition ensures that for any positive function $f\in \Ca_V(E)$ and $n\geq 1$ we have $$Q_n(f)/V\leq Q_n(V)/V \leq \Theta.$$
By (\ref{min-cond}) the integral operator $Q_n(x,dy)=q_n(x,y)~\nu(dy)$ also has a density given by $q_n(x,y)=G_{n-1}(x)p_n(x,y)$ and
 for any $r\geq r_1$ we have the local  condition 
 
\begin{equation}\label{min-cond-Q}
0<\inf_{\Va(r)^2} q_n\leq \sup_{\Va(r)^2} q_n <\infty.
\end{equation}

Consider the normalized measure valued process $\eta_n\in \Pa_V(E)$ starting at $\eta_0\in \Pa_V(E)$ defined for any $n\geq 1$ by

\begin{equation}\label{limit-ga}
\eta_{n+1}=\phi_{n+1}\left(\eta_{n}\right):=\psi_{G_n}(\eta_n)P_{n+1}\quad \mbox{\rm and}\quad
\widehat{\eta}_n:= \psi_{G_n}(\eta_n),
\end{equation}
with the updated Boltzmann-Gibbs transformations $ \psi_{G_n}$ associated with the potential function $G_n$ defined by
$$
 \psi_{G_n}(\eta_n)(dx):=\frac{1}{\eta(G_n)}~G_n(x)~\eta(dx).
$$
We readily check that the evolution semigroup $\phi_{k, n}=\phi_{k+1, n}\circ\phi_k$ associated with the flow of measure $\eta_n$ is given for any $k\leq n$ by the formula
$$
\phi_{k, n}(\eta_k)=\frac{\eta_k Q_{k, n}}{\eta_k Q_{k, n}(1)}.$$
Note that for any $\mu\in \Pa_V(E)$ we have the updating formula

\begin{equation}\label{updating_formula}
\phi_{k, n}(\eta_k)=\psi_{H^{\mu}_{k, n}}(\eta_k)\bar Q_{k,n},
\end{equation}
with the Markov operator
\begin{equation}
\bar Q_{k, n}(f) := {Q_{k, n}(f)}/{Q_{k, n}(1)}~~\text{and}~~  H_{k, n}^{\mu}(x) := \frac{Q_{k, n}(1)(x)}{\phi_{0, k}(\mu) Q_{k, n}(1)}. \label{H_n}
\end{equation}
By Theorem 4.2 in \cite{del2021stability} (see also Theorem 1 in~\cite{10.1214/12-AAP909} in the context of nonlinear filtering), for any $\mu\in \Pa_V(E)$, there exists some parameters $a,b>0$ such that for any $k\leq n$, and any $\mu_1,\mu_2\in\Pa_V(E)$ we have

\begin{equation}
\left |\left |\mu_1 \Bar Q_{k,n} - \mu_2 \Bar Q_{k, n} \right |\right |_{V} \leq ae^{-b(n-k)}||\mu_1 - \mu_2||_{V/H_{k ,n }^{\mu}} \label{contraction_h_process}.
\end{equation}
Note that for constant potential functions $G_n(x)=G_n(y)$ we have $\Bar Q_{k,n} =P_{k,n}$ and
 the above contraction estimates resume to (\ref{contraction-P}). Moreover, we stress that our focus will be on establishing exponential convergence rate without focusing on the estimation constant. Therefore, we are not concerned with whether certain constants in (3) and (10), such as \(a\), might be large due to factors like high dimensionality. Our primary concern lies solely in their existence.
 
For time homogeneous models $Q_n=Q$, Theorem 4.4 in \cite{del2021stability} ensures the existence of a leading eigen triple
$
(h, E_0, \eta_\infty)\in \left(C_V(E)\times\R_+^*\times\Pa_V(E)\right),
$
such that 
\begin{equation}\label{triplet-cond}
Q(h) = E_0~ h~~,~~\eta_\infty Q = E_0~\eta_\infty ~~\mbox{and}~~\eta_\infty(h)=1.
\end{equation}

\subsection{Schrödinger semigroups}

The objects defined in the previous subsection are core to a variety of physics problem. Indeed, consider an Hamiltonian differential operator $\Ha$ given by the formula

\begin{equation*}
\Ha := -\La + U,
\end{equation*}
where $U$ is a potential energy function from $E$ to $\R_+$, and $\La$ is a kinetic energy operator acting on a subset $\Da(\La)$ of $\Ca(E)$. 
The time dependent Schr\"odinger equation and the imaginary time version associated with the hamiltonian $\Ha$ are given, respectively, by the equations
\begin{equation*}\label{Schrodinger}
i\,\partial_t \psi_t(x)=\Ha(\psi_t)(x)\quad \mbox{\rm and}\quad
-\partial _t\varphi_t(x)=\Ha(\varphi_t)(x), 
\end{equation*}
 {with prescribed initial conditions $(\psi_0,\varphi_0)$.}
In the above display, $i\in\CC$ stands for the imaginary unit. The right-hand side equation is obtained via a formal time change by setting
$\varphi_t(x)=\psi_{-it}(x)$, and can be equivalently written in the following form
\begin{equation}\label{def-Schrod-intro-imaginary}
\partial_t \varphi_t(x)=\La(\varphi_t)(x)-U(x)\varphi_t(x), 
\end{equation}

with initial condition $\varphi_0$.

For a  twice differentiable function $\varphi_0$, the solution of (\ref{def-Schrod-intro-imaginary}) is given by the Feynman-Kac path integral formula
\begin{eqnarray}
\varphi_t(x)&=&\Qa_t(\varphi_0)(x):=\int \Qa_t(x,dy)~\varphi_0(y)\label{sol_shrodinger}\\
&=&\EE\left(\varphi_0(\Xa_t)~\exp{\left(-\int_0^t U(\Xa_s)~ds\right)}~|~\Xa_0=x\right)\nonumber.
\end{eqnarray}
In the above display, $\Xa_t$ stands for a time homogeneous stochastic process $\Xa_t$ on $E$, with generator $\La$. 
To facilitate the interpretation of the theoretical and numerical physics in the measure theoretical framework used in this article, we note that the Feynman-Kac propagator defined by the integral operator (\ref{sol_shrodinger}) is sometimes written in terms of the exponential of the Hamiltonian  operator with the exponential-type symbol
 $$
\Qa_t:=e^{-t\Ha}, $$

or in the bra-kets formalism

$$\Qa_t(\varphi_0)=\vert e^{-t\Ha}\vert\varphi_0\rangle.
$$

The integral operator $\Qa_t$ is sometimes called the Feynman-Kac propagator. 
For any $s,t\geq 0$ the integral operators $\Qa_t$ satisfy the semigroup property
\begin{align*}
    \Qa_{s+t}(x,dz)=(\Qa_s\Qa_t)(x,dy) &:=\int~\Qa_s(x,dz)~\Qa_t(z,dy) \\
    &\Longrightarrow\varphi_{s+t}=\Qa_{s}(\varphi_t).
\end{align*}

In terms of left action bra-kets, defining $\mu_{\varphi}(dx):=\varphi(x)dx$, Fubini's theorem yields
$$
\begin{array}[t]{rcl}
\displaystyle \langle \varphi\vert e^{-s\Ha}\vert \varphi_t\rangle &=&\displaystyle \int dx~ \varphi(x) ~\Qa_s(x,dy)~\varphi_t(dy)= (\mu_{\varphi}\Qa_s)(\psi_t)\\
&&\\
&=&\displaystyle  \mu_{\varphi}((\Qa_s\Qa_t)(\varphi_0)) \\
&=&\mu_{\varphi}\Qa_{s+t}(\varphi_0)  =\langle \varphi\vert e^{-(s+t)\Ha}\vert \varphi_0\rangle.
\end{array}
$$

The exponential notation is compatible with finite space models and the matrix notation of the continuous one-parameter semigroup for time homogeneous models. The bra-ket notation (a.k.a. Dirac notation) is also used to represents linear projection forms acting on Hilbert spaces associated with some reversible or some stationary measure, such as the Lebesgue measure for the harmonic oscillator.

The present article deals with different types of non necessarily stationary stochastic processes, including the free evolution process $X_t$ discussed in  \eqref{sol_shrodinger}. Apart in the reversible situation in which spectral theorems are stated on the Hilbert space associated with a reversible measure, the use of the exponential symbol or the use of the bra-kets formalism is clearly not adapted to represent different expectations with respect to different types of stochastic and non-necessarily reversible processes.

To analyze these general stochastic models, we have chosen to only use elementary and standard measure theory notation such as \eqref{action_measure}. The integral actions of a given integral operator $Q_t$ on the right for functions and on the left for measures are clearly compatible with finite space models and matrix notation. The left action $\mu \mapsto \mu \Qa_t$ maps measures into measures, while the right action $f \mapsto  \Qa_t(f)$ maps functions into functions.

\begin{equation*}
\left\{
    \begin{array}{ll}
        (\mu \Qa_t)(dy) := \int_E \mu(dx) \Qa_t(x, dy) \\ \\
        \Qa_t(f)(x) := \int_E \Qa_t(x, dy)f(y).
    \end{array}
\right.
\end{equation*}
In this context the normalized measures are defined for any $s\leq t$ by
the flow of measures
\begin{equation}\label{continu}
\mu_t(f):=\Phi_{s,t}(\mu_s)(f):={\mu_s\Qa_{t-s}(f)}/{\mu_s\Qa_{t-s}(1)},
\end{equation}
where $\mu_0$ stands for the distribution of the initial random state $\Xa_0$.
Note that
\begin{align}\label{continu-e}
\displaystyle\partial_t\mu_0\Qa_t(1) &=-
\EE\left(U(\Xa_t)\exp{\left(-\int_0^t U(\Xa_s)~ds\right)}\right) \notag \\
 &=-\mu_0\Qa_t(U) \notag \\
\displaystyle\Longrightarrow\partial_t&\log{\mu_0\Qa_t(1)}=-\mu_t(U) \notag \\
\Longrightarrow \mu_0&\Qa_t(1)=\exp{\left(-\int_0^t\mu_s(U)ds\right)}.
\end{align}

The operator $\Qa_t$ defined in \eqref{sol_shrodinger} is sometimes called a Feynman-Kac propagator. However, despite its mathematical elegance, it can rarely be solved analytically. Under some regularity conditions (cf. for instance~\cite{del2021stability}) the flow of measures $\mu_t$ converge as $t\rightarrow\infty$ to some limiting fixed point measure $\mu_{\infty}=\Phi_{s,t}(\mu_{\infty})$ (a.k.a. quasi-invariant measure). In this case, choosing $\mu_0=\mu_{\infty}$ in (\ref{continu-e}) we have
$$
\mu_{\infty}\Qa_t(1)=e^{-\lambda_0t}\quad \mbox{\rm with}\quad
\lambda_0:=\mu_{\infty}(U).
$$
Whenever it exists, the ground state $h_0$ is the leading eigen-function associated with $\lambda_0$; that is, we have
$$
\Qa_t(h_0)=e^{-\lambda_0t}~h_0.
$$
On a
given time mesh $ t_n=n\delta$, with time step $t_n-t_{n-1}=\delta>0$ we clearly have
$$
\mu_{t_n}(f)=\frac{\mu_{t_{n-1}}\Qa_{\delta}(f)}{\mu_{t_{n-1}}\Qa_{\delta}(1)}
\quad \mbox{\rm and}\quad
\Qa_{\delta}(x,dy)=G_{\delta}(x)~P_{\delta}(x,dy),$$
with the function
$$
\Ga_{\delta}(x):=\Qa_{\delta}(1)(x), $$

and the Markov transition

$$\quad \Pa_{\delta}(x,dy):=\frac{\Qa_{\delta}(x,dy)}{\Qa_{\delta}(1)(x)}.
$$
Note that
$$
\mu_{\infty}(\Ga_{\delta})=\mu_{\infty}\Qa_{\delta}(1)=e^{-\lambda_0\delta}, $$ 
as well as

$$\Qa_{\delta}(h_0)=e^{-\lambda_0\delta}~h_0
\quad \mbox{\rm and}\quad \mu_{\infty}\Qa_{\delta}=e^{-\lambda_0\delta}~\mu_{\infty}.
$$
In other words, choosing $(G,P)=(\Ga_{\delta},\Pa_{\delta})$ in the time-homogeneous Feynman-Kac model discussed in (\ref{limit-ga}) and
(\ref{triplet-cond}), we obtain the leading triple $$(\eta_{\infty},E_0,h)=(\mu_{\infty},e^{-\lambda_0\delta},h_0).$$ 
This yields the formula
$$
-\frac{1}{\delta}~\log{\mu_{\infty}(\Ga_{\delta})}=\lambda_0=\mu_{\infty}(U).
$$

Unfortunately, with the notable exception of coupled harmonic models (cf.~\cite{del2023coupled} as well as Proposition 7.1 in~\cite{arnaudon2023lyapunov}), the potential function $\Ga_{\delta}$ can rarely be evaluated and the Markov transition $\Pa_{\delta}$ cannot be sampled.
The Feynman-Kac measure $\eta_n$ introduced in (\ref{limit-ga}) can also be interpreted as the solution of a discrete-time approximation of the formula (\ref{sol_shrodinger}). Indeed, 
consider a discrete time approximation $X_{t_n}$ of the process $\Xa_{t_n}$  and let
$$
Q(x,dy)=G(x)P(x,dy),
$$
with
\begin{equation}\label{discrete-FK}
\left\{
    \begin{array}{ll}
        G(x)=\exp{(-U(x)\delta)} \\ \\
        P(x,dy)=\PP(X^{\delta}_{t_n}\in dy~|~X^{\delta}_{t_{n-1}}=x).
    \end{array}
\right. 
\end{equation}
In this situation, choosing $n=\lfloor t/\delta\rfloor$ we have 
$$
\begin{array}{l}
\displaystyle \eta_0 Q_{0,n}(f)=\mathbb{E}\left(  f(X^{\delta}_{t_n})~\prod_{0\leq k<n} G(X^{\delta}_{t_k})\right)\\
\\
\displaystyle=\mathbb{E}\left(  f(X^{\delta}_{t_n})~\exp {\left\{  -\sum_{0\leq k < n}U(X^{\delta}_{t_k}) (t_{k+1}-t_k) \right\}  }\right) \\
\simeq_{\delta\downarrow 0}\eta_0 \Qa_t(f) .
\end{array}
$$
The leading triple $(\eta_{\infty},E_0,h)$  associated with the discrete time
Feynman-Kac approximation model now depends on the time step $\delta$ and we have
$$
(\eta_{\infty},h)\simeq_{\delta\downarrow 0}(\mu_{\infty},h_0), $$

and

$$-\frac{1}{\delta}\log{\eta_{\infty}(G)}=-\frac{1}{\delta}\log{\eta_{\infty}(e^{-U\delta})}\simeq_{\delta\downarrow 0}\lambda_0=\mu_{\infty}(U).
$$

It is clearly out of the scope of this article to analyze the bias introduced by the discrete time approximation discussed above. We refer to \cite{ferre2019errorestimatesergodicproperties} for a thorough discussion of this matter.

\subsection{Diffusion Monte Carlo}\label{subsection:DMC}
The Diffusion Monte Carlo methodology relies on the fact that
the flow of measures $\eta_n$ introduced in (\ref{limit-ga}) can be interpreted as the probability distributions $\eta_n=\mbox{\rm Law}(\overline{X}_n)$ of the random states $\overline{X}_n$ of a nonlinear Markov chain $\overline{X}_n$. The choice of the Markov chain is far from unique. For instance, we have 
$$
\eta_{n+1}=\phi_{n+1}(\eta_n) = \eta_n K_{n+1,\eta_n},$$
with the local Markov transition
\begin{eqnarray*}
 K_{n+1,\eta_n}(x,dz)&:=&\left(S_{n, \eta_n}P_{n+1}\right)(x,dz)\\
 &:=&\int S_{n, \eta_n}(x,dy)~P_{n+1}(y,dz) \\
 &=&\PP\left(\overline{X}_{n+1}\in dz~|~\overline{X}_n=x\right). 
\end{eqnarray*}
In the above display, $S_{n , \eta_n}$ stands for the Markov transition
\begin{align*}
S_{n , \eta_n}(x, dy) := \epsilon_{n}(\eta_n) &G_n(x)\delta_x(dy) \\ &+ (1 - \epsilon_{n}(\eta_n) G_n(x))~\psi_{G_n}(\eta_n)(dy). 
\end{align*} 
 for some tuning parameter $ \epsilon_{n}(\eta_n)\in [0,1]$ chosen such that $ \epsilon_{n}(\eta_n)G_n(x)\in [0,1]$. For instance, for $]0,1]$-valued potential functions, we can choose
 $\epsilon_{n}(\eta_n)=0$ as well as $\epsilon_{n}(\eta_n)=1$. For more general models, we can also choose the inverse of the $\eta_n$-essential supremum of $G_n$.

 Note that the transition $\overline{X}_n\leadsto \overline{X}_{n+1}$ depends on the probability distributions $\eta_n$ of the random states $\overline{X}_n$.
 In reference with similar nonlinear Markov chain models arising in fluid mechanics, the Markov chain $\overline{X}_n$ is called a McKean interpretation of the flow of measures (\ref{limit-ga}).

The mean field particle interpretation associated with a given McKean model is defined by a 
 discrete-time system of $N$ walkers $\xi_n=\left(\xi_n^i\right)_{1\leq i\leq N}$. The system starts with $N$ independent copies of a random variable $\overline{X}_0=X_0$ with distribution $\eta_0$. Given the system   $\xi_{n}$ at some time $n\geq 0$, we sample $N$ conditionally independent walkers  $\xi_{n+1}^i$ with their respective distribution
\begin{equation*}\label{ga-def}
 K_{n+1,\eta_n^N}(\xi_{n}^i, dx)\quad \mbox{\rm with}\quad \eta^N_n:=\frac{1}{N}\sum_{1\leq i\leq N}\delta_{\xi^i_n}.
 \end{equation*}

In other words, the DMC method consists of approximating the measure $\eta_n$ by using the occupation measure $\eta^N_n$ associated with a system of $N$ walkers. The initial positions of the walkers are randomly chosen from the distribution $\eta_0$. The evolution of each walker follows then the following selection/mutation steps:

\begin{itemize}
\item Selection: We evaluate the current position $\xi_n^i$ of a walker and its potential value $G_n(\xi_n^i)$. With probability $\left(1 - \epsilon_{n}(\eta_n^N)G_n(\xi_n^i)\right)$, $\xi_n^i$ is killed and instantly replaced by another walker say $\widehat{\xi}_n^i=\xi^j_n$ with a probability proportional to $G_n(\xi^j_n)$ and $j\in\{1,\ldots,N\}$; otherwise we keep it and set $\widehat{\xi}^i_n:=\xi_n^i$.
\item Mutation: We move the selected walker $\widehat{\xi}_n^i=x$ to a new location $\xi_{n+1}^i=y$ using the transition kernel $P_{n+1}(x,dy)$.
\end{itemize}

 The selection transition associated with the choice  $\epsilon_{n}(\eta_n^N)=0$ coincides with the so-called proportional selection/reconfiguration. Note that the walker with the highest potential value is always selected when  $\epsilon_{n}(\eta_n^N)$ is the inverse of the $\eta^N_n$-essential supremum of $G_n$. 
 
 For $[0,1]$-valued potential functions $G_n$ we can also choose  $\epsilon_{n}(\eta_n^N)=1$. In this situation, particle are killed at a geometric clock that depends on the potential function. 
 In the context of the discrete time approximating Feynman-Kac models discussed in (\ref{discrete-FK}), when the time step $\delta$ tends to $0$, these geometric  killing-rates converges  to an exponential killing rate. The limiting DMC scheme in continuous time consists of a system of $N$ walkers. Between killing times, walkers explore the space  with a free evolution with generator $\La$; at rate $U$ the walkers are killed and instantly another walker in the pool duplicates (see for instance~\cite{del2003particle,DelMoral2000,MORAL2000193,10.1214/20-EJP546,doi:10.1137/050640667}). {\em We underline that without the geometric killing rates, the variance of the proportional reconfiguration associated with the choice  $\epsilon_{n}(\eta_n^N)=0$ blows up when the time step tends to $0$.}

We expect that the occupation measures of the system approximate the solution of the measure-valued process (\ref{limit-ga}); that is, in a sense to be given, for any time horizon $n\geq 0$ we have 
$$
\eta^N_{n}=\frac{1}{N}\sum_{i=1}^N\delta_{\xi_n^i}\underset{N\rightarrow\infty}{\longrightarrow}\eta_n,
$$
as well as
$$
\displaystyle\psi_{G_n}\left(\eta^N_{n}\right)=\sum_{i=1}^N\frac{G_n(\xi^i_n)}{\sum_{1\leq j\leq N}G_n(\xi^j_n)}~\delta_{\xi_n^i}\underset{N\rightarrow\infty}{\longrightarrow}\psi_{G_n}(\eta_n)=\widehat{\eta}_n.
$$
 
 The Lyapunov drift condition (\ref{lyap-cond}) combined with the local minorization condition (\ref{min-cond}) ensures that the free evolution of the walkers is stable, in the sense that it forgets exponentially fast its initial condition (see (\ref{contraction-P})). This condition is not satisfied for linear gaussian processes with an unstable drift matrix. In this context, as shown in Proposition~\ref{thm:div_oscillo} and Proposition~\ref{prop-fi} (see also Section~\ref{subsection:hamo_1_d})  {\em the DMC method diverges for any fixed number of walkers,  even if the asymptotic variance of the Central Limit Theorem stated in Lemma~\ref{lmm:TCL } is uniformly bounded w.r.t. the time horizon}.
 
Note that the stability regions and the Lyapunov functions are connected to the potential function by (\ref{GQV}). Importance sampling techniques and twisted guiding waves can be used to design stable-like free evolutions. For instance, using the decomposition
 $$
\widehat{Q}_{n+1}(f)(x)=\widehat{G}_{n}~\widehat{P}_{n+1}(f) ,
$$
with the potential function
$$\widehat{G}_{n}:=P_{n+1}(G_{n+1})\quad
\mbox{\rm and}\quad \widehat{P}_{n+1}(f) :=\frac{P_{n+1}(G_{n+1}f)}{P_{n+1}(G_{n+1})},
$$
we readily check the evolution equation

\begin{equation}
  \widehat{\eta}_{n+1}= \widehat{\phi}_{n+1}(\widehat{\eta}_{n}):=\psi_{\widehat{G}_{n}}(\widehat{\eta}_{n})\widehat{P}_{n+1}. \label{idem-eq-w}
\end{equation}

 This shows that the updated measures evolve  as in (\ref{limit-ga}) by replacing $(G_n,P_n)$ by $(\widehat{G}_n,\widehat{P}_n)$. The DMC associated with these objects is a  genetic-type Monte Carlo sampler with selection fitness  functions $\widehat{G}_n$ and mutation transitions $\widehat{P}_n$. From the mathematical viewpoint, this model coincides with the one discussed in (\ref{limit-ga}). Nevertheless, the free evolution of the walkers associated with $(\widehat{G}_n,\widehat{P}_n)$ is now driven by the potential function. This local conditioning importance sampling strategy if often used to turn infinite energy absorbing wells (a.k.a. hard obstacles) into soft ones~\cite{del2004particle,moral2018note}.
 
 More generally, for any given time mesh $k_n\leq k_{n+1}$ we have
 $$
  \widehat{Q}_{k_n, k_{n+1}}(f)= \widehat{G}_{k_n}~ \widehat{P}_{k_n, k_{n+1}}(f),
$$
with the potential function
$$
\widehat{G}_{k_n}:=  \widehat{Q}_{k_n, k_{n+1}}(1)\quad
\mbox{\rm and}\quad \widehat{P}_{k_n, k_{n+1}}(f) :=\frac{  \widehat{Q}_{k_n, k_{n+1}}(f)}{  \widehat{Q}_{k_n, k_{n+1}}(1)}.
 $$
 This yields the formula
 $$
 \widehat{\eta}_{k_{n+1}}(f)=\widehat{\phi}_{k_n,k_{n+1}}(\widehat{\eta}_{k_n})=
 \psi_{\widehat{G}_{k_n}}(\widehat{\eta}_{k_n})\widehat{P}_{k_n, k_{n+1}}.
 $$
  This shows that the updated measures $\widehat{\eta}_{k_n}$ evolve  as in (\ref{limit-ga}) by replacing $(G_n,P_n)$ by $(\widehat{G}_{k_n},\widehat{P}_{k_n, k_{n+1}})$. The DMC associated with these objects is a  genetic-type Monte Carlo sampler with selection fitness  functions $\widehat{G}_{k_n}$ and mutation transitions $\widehat{P}_{k_n, k_{n+1}}$.

As shown in Section~\ref{subsection:importance sampling} (see also Corollary~\ref{cor:optimal_filter}) in the context of coupled harmonic oscillators
there exists a time mesh for which the mutation transitions $\widehat{P}_{k_n, k_{n+1}}$ and the potential functions  $\widehat{G}_{k_n}$ satisfy the required stability properties. 
For more general models, these objects do not have an analytic form. In this context, we can use the
unbiased Monte Carlo methodologies  discussed in 
~\cite{del1998measure}, see also Section 2.3.2 in~\cite{DelMoral2000}, and Section 11.5 in~\cite{bookFK}.
 
 \section{Statement of the main results}\label{section:main_result}
\subsection{Some
regularity conditions}\label{subsection:regularity_cond}

For $f\in\Ca_b(E)$, time-uniform $L_p$-convergence of the error made by the DMC method in estimating $\eta_n(f)$ have been obtained (see for example \cite{DelMoral2000,del2001stability}, as well as Chapter 4 in~\cite{bookFK}, Chapter 12 in~\cite{del2013mean} and the more recent article~\cite{del2021stability}) under the strong mixing assumption that there exist $\epsilon_P\in\R_+^{*}$ and $\epsilon_G\in\R_+^{*}$ such that

$$ 
\left\{
    \begin{array}{ll}
        \forall (x_1, x_2, n)\in E^2\times \N,~~G_{n}(x_1) \ge \epsilon_G ~G_{n}(x_2) \\ \\
        P_{n}(x_1, dy) \ge \epsilon_P ~P_{n}(x_2,dy).
    \end{array}
\right. 
$$
 
A significant consequence of this assumption is a time uniform bound on the potential function defined for some $\gamma \in \Pa_V(E)$ and $k\in\N$ by

$$ G^\gamma_{k,k+n} :x\in\E \mapsto  G^\gamma_{k,k+n}(x):={1}/{H^\gamma_{k, k+n}(x)}.$$
Unfortunately, these uniform minorization
and majorization conditions are rarely satisfied when $E$ is non-compact.

In order to guarantee a time-uniform $L_p$-convergence for more general models including coupled harmonic oscillators, our framework requires to estimate uniformly in time the inverse moments of $\eta_{q}^N(H_{q, n}^\gamma)$. To do this, we first assume that there exists $\gamma \in \Pa_V(E)$ such that

\begin{equation}
	\sup_{n\in\N} ~ \phi_{0, n}(\gamma) \left (G_n\right ) < +\infty.\label{control_phi_G}
\end{equation}

For time-homogeneous models, without any further conditions on the potential, condition \eqref{control_phi_G} is easily checked with $\gamma = \eta_\infty$. For this scenario, we will then consider in the rest that $\gamma = \eta_\infty$. Moreover, this hypothesis trivially holds if the functions $G_n$ are bounded by some constant independent of $n$.

We assume that there exists $W\in\Ca_V(E)$, $\alpha \in ]0, 1]$ and a $Q$-Lyapunov function $\bar W\in \Ca_V(E)$ such that

\begin{equation}
Q_n(W) \ge \Pbar \times W ~~\text{and}~~ W^{-\alpha} \leq \Bar W ,\label{Q(W) > W general}
\end{equation}

where $$\Pbar := \sup\limits_{(n, x)\in(\N^*\times E)} P_n(G_n)(x). $$

For time homogeneous models, this condition can be relaxed into the following

\begin{equation}
Q(W) \ge \min \left \{\Pbar, E_0 \right \} \times W ~~\text{and}~~ W^{-\alpha} \leq \Bar W .\label{Q(W) > W homogeneous}
\end{equation}

Note that, for time-homogeneous models, the set of functions $W\in\Pa_V(E)$ such that  $Q(W) \ge E_0 W$ is non-empty as it contains at least the ground state $h$. It is also worth noting that it is not necessary to know the exact value of $\Pbar$ nor the one of  $E_0$ in order to prove that \eqref{Q(W) > W general} or \eqref{Q(W) > W homogeneous} hold. Indeed, if one of these constants is less than some $C\in\R\cup\{+\infty\}$, then it is sufficient to prove that for any $c < C$, there exist $W_c\in\Ca_V(E)$ and $\alpha_c\in(0,1]$ such that

$$Q(W_c) \ge c~ W_c ~~\text{and}~~ W_c^{-\alpha_c} \leq V. $$

Finally, we assume that there exists a $Q$-Lyapunov function $\Bar V\in\Ca(E)$ and $\lambda\in 2\N$ such that
 
\begin{equation}\label{lambda-reff}
V^{\lambda} \leq \Bar V.
\end{equation}

Without further mention, we assume that $V, \bar V$ and $\bar W$ are integrable with respect to $\eta_0$, i.e., $\eta_0\in\Pa_V(E)\cap \Pa_{\Bar W}(E)\cap \Pa_{\Bar V}(E)$. Under conditions \eqref{control_phi_G} and \eqref{Q(W) > W general} or \eqref{Q(W) > W homogeneous}, it is then possible to obtain a time uniform bound on the random potential function $G_{k, n}^\gamma$.

\begin{lmm}\label{lmm:control_H_l}
 Let $\gamma$ and $(\eta, \mu)$ be defined as in~\eqref{control_phi_G}. There exists $\Bar \beta \in \R^*_+$ such that for any $\beta \leq \bar{\beta}$ we have

\begin{equation}
\left\{
    \begin{array}{ll}
        \sup\limits_{\substack{(q,n, N)\in\N^3 \\ q\leq n}} \E\left [\phi_{q}(\eta_{q-1}^N)(H_{q, n}^\gamma)^{-\beta}\right ] < +\infty \\ \\
        \sup\limits_{\substack{(q,n, N)\in\N^3 \\ q\leq n}} \E\left [\eta_{q}^N(H_{q, n}^\gamma)^{-\beta}\right ] < +\infty.
    \end{array} 
\right.
\end{equation}
\end{lmm}

The proof of this pivotal Lemma is postponed to the appendix.

\subsection{A time-uniform convergence Theorem}

The main goal of this paper is to establish that in the context we described, the $L_p$-norm of the error made by the DMC method in approximating the Feynman-Kac measures $\eta_n$ remains bounded in time and converges to zero as the number of particles increases. Our main result can be stated as follows

\begin{thm} \label{thm:optimal_predictor_unif_bound}
For any  $p\in\N^*$, there exists $c\in\R$ and $\beta \in (0, 1]$ such that for any $f\in\Ca_{V^\frac{\lambda}{4p}}(E)$  we have

\begin{equation*}
\sup_{n\in\N}\E\left (|\eta_n(f) - \eta_n^N(f)|^p\right )^{\frac{1}{p}} \leq cN^{-\frac{\beta}{2}}.
\end{equation*} 
\end{thm}

The proof of this theorem is provided in subsection \ref{subsection:proof_th_1}. Note that, due to the \(\beta\) factor, the convergence rate obtained here may be slower than other existing convergence results cited previously. The presence of this factor arises from technical reasons in the proof of Theorem~\ref{thm:optimal_predictor_unif_bound} and is directly related to the choice of the Lyapunov functions \(V\) and \(W\). It could be interesting to conduct a numerical study to determine whether there are cases exhibiting a convergence rate proportional to \(N^{-\frac{\beta}{2}}\) with \(\beta\) lower than one.

For time-homogeneous models, a direct consequence of Theorem \ref{thm:optimal_predictor_unif_bound} is a control over the estimation of the limiting quasi-invariant measure $\eta_\infty$.

\begin{cor}\label{cor:limit_bound}
For any  $p\in\N^*$, there exists $(c_1, c_2, \omega)\in\R_+^{*^3}$ and $\beta \in (0, 1]$ such that for any $f\in\Ca_{V^\frac{\lambda}{4p}}(E)$ we have

$$ \E(|\eta_n^N(f) - \eta_\infty(f)|^p)^\frac{1}{p} \leq c_1{N^{-\beta/2}} +  c_2e^{-\omega n}.$$

Moreover, there exists $c\in\R$ such that, with $a = \frac{1}{\omega}\ln(c_2/c_1)$ and $b = \frac{\beta}{2\omega}$, have

\begin{equation*}
\sup_{n\ge a + b\ln(N)}\E(|\eta_n^N(f) - \eta_\infty(f)|^p)^\frac{1}{p} \leq c N^{-\frac{\beta}{2}}.
\end{equation*}
\end{cor}

The proof the above Corollary is provided in subsection \ref{subsection:grd_state_approx}
.

Assuming that there exists a measure $\mu\in\Pa(E)$ that is reversible for $P$, it becomes possible to obtain a re-normalized weak form of the ground state $h$ and its associated eigenvalue from the limit measures $\eta_\infty$ and $\widehat\eta_\infty$ of $\eta_n$ and $\widehat\eta_n$.

Indeed, referring to Section~9.5.5 in \cite{delmoral:halartilce1} (see also~\cite{del2003particle} as well as Section 12.4 in~\cite{bookFK}), we have
\begin{equation*}
\eta_\infty(G) = E_0, 
\end{equation*}

and

\begin{equation*}
\eta_\infty(f) = \frac{\mu(P(h)f)}{\mu(P(h))}=\frac{\mu(hP(f))}{\mu(h)}=\psi_h(\mu)P(f).\label{limit_objects}
\end{equation*}
Note that
$$
\eta_\infty(f) = \frac{\mu(Q(h)f/G)}{\mu(h)}=E_0~\psi_h(\mu)(f/G).
$$
In the reverse angle, we have the updated limiting measures
$$
\widehat \eta_\infty(f) = \psi_{G}(\eta_\infty)(f)=\psi_h(\mu)(f)= \frac{\mu(hf)}{\mu(h)} \Longrightarrow
 \eta_\infty=\widehat \eta_\infty P.$$

The existence of a reversible measure is actually not required to express the ground-state energy using the limit measure; indeed, we always have

\begin{equation*}
\eta_\infty(h) = \phi(\eta_\infty)(h) = \frac{\eta_\infty(GP(h))}{\eta_\infty(G)} = E_0~\frac{\eta_\infty(h)}{\eta_\infty(G)}~\Rightarrow~ \eta_\infty(G) = E_0.
\end{equation*}

Those equalities, combined with the convergence stated in Corollary \ref{cor:limit_bound}, guarantee the efficiency of the DMC method for approximating the ground-state energy and wave function of quantum systems.

\begin{cor}\label{cor:control_limit}
Let $p\in\N^*$, and assume that $G\in\Ca_{V^\frac{\lambda}{4p}}(E)$. There exists $(a, b, c)\in\R_+^{*^3}$ and $\beta \in (0, 1]$ such that for any $f\in\Ca_{V^\frac{\lambda}{4p}}(E)$ we have

\begin{equation*}
\left\{
    \begin{array}{ll}
        \sup\limits_{n\ge a + b\ln(N)}\E\left (|E_0 - \eta_n^N(G)|^p\right )^{\frac{1}{p}} \leq cN^{-\frac{\beta}{2}} \\
        \sup\limits_{n\ge a + b\ln(N)}\E\left (\left |\psi_h(\mu)(P(f)) - \eta_n^N(f)\right |^p\right )^{\frac{1}{p}} \leq cN^{-\frac{\beta}{2}} \\
        \sup\limits_{n\ge a + b\ln(N)}\E\left (\left |E_0~\psi_h(\mu)(f/G) - \eta_{n}^N(f)\right |^p\right )^{\frac{1}{p}} \leq cN^{-\frac{\beta}{2}}  .
    \end{array}
\right. 
\end{equation*} 
\end{cor}

\subsection{Coupled harmonic oscillators}\label{subsection:harmonic_oscillator}

To illustrate the practical applications of Theorem \ref{thm:optimal_predictor_unif_bound}, we carry out an in-depth study of the generalized coupled harmonic oscillator~\cite{del2023coupled}. Given $d\ge 1$, consider the flow of measures $\mu_t$ on $\R^d$ defined by \eqref{continu} with the the kinetic energy operator $\La$ and the potential function $U$ defined for some $d\times d$ real matrices $(C, D, F)$ and $x\in\R^d$ by:

$$ U(x) := \exp\left (-\frac{1}{2}x^\prime Fx\right )$$

and

$$\La := \sum\limits_{1\leq k,l \leq d}C_{k,l}~x_l~\partial_{x_k} + \sum\limits_{1\leq k,l \leq d} D_{k,l}~\partial_{x_k,x_l}, $$

with $x^\prime$ standing for the transposition of vector $x$ (or matrix, where appropriate) .

We define the normalized Markov integral operator $\bar{\Qa}_t$, associated with \(\Qa_t\) from equation \eqref{sol_shrodinger}, using the ratio formula:

$$\bar{\Qa}_t(f)(x) := \Qa_t(f)(x)/\Qa_t(1)(x)$$

We know from ~\cite{del2023coupled} that, for $t>0$, there exists two positive-definite matrices $B$ and $S$, a matrix $A$ and $C\in\R_+^*$ such that: 
$$\left\{\begin{array}{ll}
        \bar{\Qa}_t(f)(x) \sim \mathcal{N}(Ax, B) \\
        Qa_t(1)(x) = C\exp\left (-\frac{1}{2}x^{\prime} S x\right ) 
    \end{array}.  \right.
$$

For any $n\in \mathbb{N}$, the sequence of measure $\gamma_{tn}$ is then equal to the flow of measures $(\eta_n)_{n\in\N}$ starting from $\gamma_0$ and evolving according to \eqref{limit-ga} with the the transitions kernel $P_{A, B}$ and the positive function  $G_{S}$ defined by

\begin{equation}\label{def_P_A_B}
\left\{
    \begin{array}{ll}
        P_{A, B}(x, dy) =  \frac {e^{-{\frac {1}{2}}(Ax-y )^{\prime }B ^{-1}(Ax-y )}}{(2\pi )^{k/2}\left|B \right|^{1/2}}dy \\ \\
        G_{S}(x):= \exp\left (-\frac{x^\prime Sx}2\right )
    \end{array}.
\right.
\end{equation}

For the harmonic oscillator, studying the discrete-time models instead of the continuous-time models involves then no approximation. Therefore, we will focus on discrete-time models in the following analysis.

Note that, strictly speaking, in physical terms, the harmonic oscillator framework implies that $\eta_0$ is a normal distribution and that $A$ is the null matrix. However, for the purposes of this study, within the coupled harmonic oscillator framework, we assume that the initial measure is not restricted to a normal distribution and that  $A$ can be any real matrix, not necessarily the null matrix. Additionally, the matrix $A$ may not be symmetric nor a stable (Hurwitz) matrix. Furthermore, the transition $P_{A, B}$ is not necessarily reversible, unless $AB=BA^{\prime }$.

In order to state a convergence result on the error of the DMC method that is not restricted to the harmonic oscillator, we now consider a family of Feller Markov transitions $(P_n)_{n\in\N}$ with  positive densities $p_n$, an initial distribution $\eta_0\in\Pa(E)$ and a family of positive functions $(G_n)_{n\in\N}\in\Ca_0(\R^d)^\N$ bounded by $1$.

Let $(A, \mathfrak{A})$ be real matrices and $(B, \mathfrak{B}, S)$ be real positive definite matrices. Denote by $p_{A, B}$ the density of $P_{A, B}$ and by $E_{A, B, S}$ the ground-state energy associated with the operators $P_{A, B}$ and $G_S$. We assume then the following

\begin{itemize}
\item There exists $c_1 \in \R$ such that for any $(n , x, y)\in\N\times\R^{d}\times\R^{d} $ we have
\begin{equation}\label{upper_bound_density}
p_n(x, y) \leq c_1~ p_{A, B}(x, y). 
\end{equation} 
\item For any $(n, x, y)\in\N\times \R^{d} \times \R^d$ we have
\begin{equation}\label{lower_control} G_n(x)~p_{n+1}(x, y) \ge \frac{\Pbar }{E_{\mathfrak{A}, \mathfrak{B}, S}}\times G_{S}(x)~p_{\mathfrak{A}, \mathfrak{B}}(x, y).
\end{equation}
\item There exists a compact $K \subset \R^d$ such that for any $ (n ,x, y)\in\N\times\R^{d}\times \R^d~$ we have \begin{equation} \label{variation_G_GS}
[G_S(x) - G_S(y)][G^{-1}_n(x)\mathbb{1}_{\R^d\setminus K}(x) - G^{-1}_n(y)\mathbb{1}_{\R^d\setminus K}(y)] \ge 0.
\end{equation}
\item There exists $c_3 \in \R_+^*$ such that $G^{-c_3}_S$ is integrable with respect to $\eta_0$.
\end{itemize}

If $P_n$ and $G_n$ are time-independent, it is possible to replace $\Pbar$ in \eqref{lower_control} by the ground state energy associated with $P$ and $G$.

These conditions hold trivially for the coupled harmonic oscillator, i.e if $P_n = P_{A , B}$, $G_n = G_S$ and $\eta_0$ is a normal distribution.

In this context, Theorem \ref{thm:optimal_predictor_unif_bound} leads to a simple sufficient matrix condition that guarantees the uniform convergence of the DMC method. The proof of the following corollary can be found in subsection \ref{subsection:hamo_mutli_d}.

\begin{cor}\label{cor:conv_osci}
Consider a family of Feller Markov transitions $(P_n)_{n\in\N}$  and a family of positive functions $(G_n)_{n\in\N}\in\Ca_0(\R^d)^\N$  bounded by $1$, verifying \eqref{upper_bound_density}, \eqref{lower_control} and \eqref{variation_G_GS} with
$A^\prime SA < S$. Then, for any $p \in \N^*$, there exist $(\beta, \alpha, c) \in \left(0, 1\right] \times \R_+^{*^2}$ such that for any function $f \in \Ca_V(\R^d)$ we have

\begin{equation*}
\sup_{n\in\N}\E\left (|\eta_n(f) - \eta_n^N(f)|^p\right )^{\frac{1}{p}} \leq cN^{-\frac{\beta}{2}}, \end{equation*} 
with
\begin{equation*}
V: x\in\R^d \mapsto \exp\left (\frac{\alpha}{2}~x^TSx\right ).
\end{equation*} 
\end{cor}

Returning to the coupled harmonic oscillator, we consider $P = P_{A, B}$ and $G = G_S$. Assuming the same conditions on $\eta_0$ as in Corollary \ref{cor:conv_osci}, we aim to design an algorithm that estimates $\eta_n$ and such that its error converges to $0$ for any value of $A$, $B$ and $S$. To satisfy the convergence condition of Corollary \ref{cor:conv_osci}, we propose a modification of the DMC method by introducing a change in the transition kernel and the selection function.  We define, for all $k \geq 1$, the functions $G^{(k)} \in \mathcal{C}_V(E)$ and the Markov kernels $P^{(k)}$ on $E$, such that for all $f \in \mathcal{C}_V(E)$, we have:

\begin{equation}
Q_{0, k}(f) =   G^{(k)} P^{(k)}(f)\quad\mbox{\rm with}\quad   \left\{
    \begin{array}{ll}
         G^{(k)} :=  Q_{0, k}(1) \\ \\
        P^{(k)}(f) := { Q_{0, k}(f)}/{ Q_{0, k}(1)} .
    \end{array}
\right. \label{kernel_analogy}
\end{equation} 
We use the convention $  G^{(0)}=1$ and $P^{(0)}=Id$ for $k=0$.

For $k\in\N^*$, consider a system of walkers $\xi_n^{(k)}=\left(\xi_n^{(k), i}\right)_{1\leq i\leq N}$ associated to the DMC method with initial distribution $\psi_G(\eta_0)$, transitions ${P}^{(k)}$ and selection function ${G}^{(k)}$ as well as the empirical measures

\begin{equation}
 \eta^{(k), N}_{n}:=\frac{1}{N}\sum_{1\leq i\leq N}\delta_{ \xi_n^{(k), i}}. \label{approx_hat}
\end{equation}

This system of walkers offers an approximation of the measures ${\eta}_n$ for any $n\in k\N$. This type of change in the approximation, based on an importance sampling transformation is analogous to using a guiding waves function to direct the Monte Carlo moves. Without any additional condition, Theorem \ref{thm:optimal_predictor_unif_bound} ensures the uniform convergence of the model. The details of the proof can be found in Subsection \ref{subsection:hamo_mutli_d}. 
 
\begin{cor}\label{cor:optimal_filter}
Let $p \in \N^*$, there exists $\bar k \in \N$ such that for any $k\ge\bar k$, there exist $(\beta, \alpha, c) \in \left(0, 1\right] \times \R_+^{*^2}$ satisfying  for any  $f \in \Ca_V(\R)$

\begin{equation*}
\sup_{n\in\N}\E\left (|\widehat \eta_{nk}(f) - \eta_n^{(k), N}(f)|^p\right )^{\frac{1}{p}} \leq cN^{-\frac{\beta}{2}}, \end{equation*}

with 
\begin{equation*}
V: x\in\R^d \mapsto \exp\left (\frac{\alpha}{2}x^TSx\right ).
\end{equation*} 
\end{cor}

Although the method does not provide an approximation for every time step, several strategies can be used to fill the gaps left by the approximation. A simple approach, though more computationally intensive,  is to run independent systems of walkers for each time step in the interval $\llbracket 0, k-1 \rrbracket$. This method not only fills in the gaps, but also maintains the convergence property.

So far, our focus has been on approximating the flow of measures $\eta_n$. However, as shown in \eqref{idem-eq-w}, the flow of measures $\widehat{\eta}_n$ is defined similarly to $\eta_n$ but with $(G,P)$ replaced by $(\widehat{G},\widehat{P})$. In the context of the coupled harmonic oscillator discussed earlier,  denoting by $\mathcal{N}(A , B)$ the Gaussian distribution centered at $A$ with co-variance matrix $B$,  we have:

\begin{equation}
\left\{
    \begin{array}{ll}
        \delta_x\widehat{P}=\Na(\widehat{A}x,\widehat{B}) \\ \\
        \widehat{G}(x)\propto \exp{\left(-\frac{1}{2}~x^{\prime }\widehat{S}~x\right)}
    \end{array}
\right. 
\end{equation}
with
\begin{equation}
\left\{
    \begin{array}{ll}
\widehat{B}:=(B^{-1}+S)^{-1}\leq B \\ \widehat{A}:=\widehat{B}~B^{-1}A \\
\widehat{S}:=A^\prime \left(S-S(B^{-1}+S)^{-1}S\right)A\geq 0.
    \end{array}
\right. \label{hat P and G}
\end{equation}

This shows that, within this framework, the convergence condition stated in Corollary~\ref{cor:conv_osci} can be extended to the updated measures $(\widehat \eta_n)_{n\in\N}$ by replacing $A$ and $S$ with $\widehat{A}$ and $\widehat{S}$. Similarly, Corollary~\ref{cor:optimal_filter} can be extended in the same way.

Note that all corollaries in this subsection can be extended to a control on the estimation of the limit measures, ground state, and eigenvalue using the same approach as presented in Corollary \ref{cor:control_limit}.

Our study concludes with a focus on the divergence of the DMC method when approximating the one-dimensional harmonic oscillator. This confirms that the stability condition stated in Corollary \ref{cor:conv_osci} is necessary and that, in some cases, the set of assumptions presented can closely approximate a sufficient and necessary condition. Additionally, it emphasizes the significance of the importance sampling method introduced in the previous corollary. Specifically, in the one-dimensional context, the sufficient condition for uniform convergence of the DMC method is expressed by $A^2 < 1$. Proposition \ref{thm:div_oscillo} establishes the divergence of the error made by the DMC method when $A^2 > 1$, leaving open only the case $A=1$.

\begin{prop}\label{thm:div_oscillo}
Assume that $A^2>1$ and $P_0>0$. For any $p\in\N*$ we have

\begin{equation*}
\sup_{n\in\N} \E\left (|\eta_n(I) - \eta_n^N(I)|^p\right )^{\frac{1}{p}} = +\infty.
\end{equation*}
\end{prop}

We complement the divergence result with an observation on the asymptotic convergence. Our next proposition proves that, regardless of the value of $A$, the univariate harmonic oscillator can always be uniformly controlled over time as the number of walkers approaches infinity. This underscores the importance of achieving uniform control for a finite number of walkers, rather than focusing solely on the asymptotic scenario.

\begin{prop}
Consider the univariate harmonic oscillator.
For any $n\ge 1$, we have the following convergence in law as N tends toward $+\infty$

$$
\sqrt{N}[\eta_n(I) - \eta_n^N(I)] \xrightarrow[N \rightarrow \infty]{\mathcal{L}} \mathcal{N}(0 , \sigma_n^2) ,
$$

with $\sigma^2_n$ the asymptotic variance such that

$$
\sup_{n\in\N} \sigma^2_n < \infty. 
$$

\end{prop}

The proof of those propositions can be found in Subsection \ref{subsection:hamo_1_d}.

\section{Stochastic interpolation} \label{section:proof_bound}

\subsection{Time varying semigroups}\label{subsection:proof_th_1}

In this subsection, we focus on proving Theorem \ref{thm:optimal_predictor_unif_bound}. To take advantage of the conditional independence of the walkers, we structure our approach around the following decomposition of the difference between the Feynman-Kac measure and its empirical approximation, using the convention $\eta_{-1}^N = \eta_{0}$. Following~\cite{DelMoral2000,del2001stability}, we use the following stochastic interpolation formula

\begin{equation}
\eta^N_n - \eta_n = \sum\limits_{q=0}^n [\phi_{q, n}(\eta_q^N) - \phi_{q, n}(\phi_{q}(\eta_{q-1}^N))]. \label{ont_step_decomp}
\end{equation}

Each term on the right-hand side represents the error that occurs when using the DMC approximation instead of the real propagator for a single extra time step. Combining the uniform bound given in Lemma \ref{lmm:control_H_l} with the contraction property \eqref{contraction_h_process}, the following Lemma establishes an exponentially decreasing control for these local errors.
 
\begin{lmm}\label{lmm:bound_one_step}
For any  $p\in\N^*$, there exists $(c,\rho, \beta)\in\R_+^{*^2}\times(0,1]$ such that for any function $f\in\Ca_{V^{\frac{\lambda}{4p}}}(E)$  and any $(N, q, n)\in\N^3$ with $q \leq n$ we have

\begin{equation}
\E\left[\left |[\phi_{q, n}(\eta_q^N) - \phi_{q, n}(\phi_{q}(\eta_{q-1}^N))](f) \right |^p \right ]^{\frac{1}{p}} \leq ce^{-(n-q)\rho}N^{-\frac{\beta}{2}}. \label{time_n}
\end{equation}
\end{lmm}

\textbf{Proof}: \\

Let $(\eta, \mu)\in\Pa(E)$ and let $\gamma\in\Pa_V(E)$ be defined as in \eqref{control_phi_G}. Consider  $H_{q, n} := H^\gamma_{q, n}$ as defined in \eqref{H_n}. Applying the updating formula \eqref{updating_formula}, we obtain
 
 \begin{align*}
&\phi_{q, n}(\eta)(f) - \phi_{q, n}(\mu)(f) \\ 
&= (\psi_{H_{q, n}}(\eta)\Bar Q_{q, n} - \psi_{H_{q, n}}(\mu)\Bar Q_{q, n})(f) \notag \\
&= \frac{1}{\eta(H_{q, n})}(\eta - \mu)\left (H_{q, n}\Bar Q_{q, n}[f - \psi_{H_{q, n}}(\mu)\Bar Q_{q, n}(f)]\right ).  \label{difference_psi_M}
\end{align*}
This yields the formula
\begin{equation}\label{kr-m}
\phi_{q, n}(\eta)(f) - \phi_{q, n}(\mu)(f)=\frac{1}{\eta(H_{q, n})}(\eta - \mu)(F^\mu_{q, n}),
\end{equation}
with the function
$$
F_{q, n}^\mu(x) := H_{q, n}(x)\int_E\psi_{H_{q, n}}(\mu)(dy)[\Bar Q_{q, n}(f)(x) - \Bar Q_{q, n}(f)(y)] .$$

Then, applying Hölder's inequality for any $\beta \in [0, 1)$, we obtain the estimate

\begin{align*}
&\E\left (\left |\phi_{q, n}(\eta_q^N)(f) - \phi_{q, n}(\phi_{q}(\eta_{q-1}^N))(f)\right |^p\right )^{\frac{1}{p}}  \notag \\
\leq ~ &\E\left (\eta_q^N(H_{q, n})^{-2p}\left |(\eta_q^N - \phi_{q}(\eta_{q-1}^N))(F_{q, n}^N)\right |^{2p(1-\beta)}\right )^{\frac{1}{2p}} \\ &\times\E\left (\left |(\eta_q^N - \phi_{q}(\eta_{q-1}^N))(F_{q, n}^N)\right |^{2\beta p}\right )^{\frac{1}{2p}}, \notag 
\end{align*} 
with $F^N_{q, n} := F_{q, n}^{\phi_{q}(\eta_{q-1}^N)}$. 
Using  (\ref{kr-m}) this yields the estimate
$$
\begin{array}{l}
\E\left (\left |\phi_{q, n}(\eta_q^N)(f) - \phi_{q, n}(\phi_{q}(\eta_{q-1}^N))(f)\right |^p\right )^{\frac{1}{p}}  \notag \\
\leq ~~\E\left (\left |(\eta_q^N - \phi_{q}(\eta_{q-1}^N))(F_{q, n}^N)\right |^{2\beta p}\right )^{\frac{1}{2p}} \times
\\
\E\left (\eta_q^N(H_{q, n})^{-2\beta p}\left |\phi_{q, n}(\eta_q^N)(f) - \phi_{q, n}(\phi_{q}(\eta_{q-1}^N))(f)\right |^{2p(1-\beta)}\right )^{\frac{1}{2p}}.
\end{array}
$$
Recalling that $f\in\Ca_{V^{\frac{\lambda}{4p}}}(E)$, this implies that
$$
\begin{array}{l}
\E\left (\left |\phi_{q, n}(\eta_q^N)(f) - \phi_{q, n}(\phi_{q}(\eta_{q-1}^N))(f)\right |^p\right )^{\frac{1}{p}} \\
\\
\leq 
\left [\E\left (\phi_{q, n}(\eta_q^N)(V^\lambda)\right )^{\frac{1}{4p}}+ \E\left (\phi_{q, n}(\phi_{q}(\eta_{q-1}^N))(V^\lambda)\right )^{\frac{1}{4p}}\right ]\\
\\
\times \E\left (\eta_q^N(H_{q, n})^{-4\beta p}\right )^{\frac{1}{4p}}\E\left (\left |(\eta_q^N - \phi_{q}(\eta_{q-1}^N))(F_{q, n}^N)\right |^{2\beta p}\right )^{\frac{1}{2p}}.
\end{array}
$$
From Lemmas \ref{lmm:control_H_l} and \ref{lmm:unif_bound_lyapu}, we deduce that, to conclude, it is enough to prove that, for some constant $c\in\R_+^{*}$ independent of $n$, $q$ and $N$, we have

\begin{equation}
\E\left [\left |\left [\eta_q^N - \phi_{q}(\eta_{q-1}^N)\right ]\left (F^N_{q, n}\right )\right |^{2\beta p}\right ]^{\frac{1}{2p}}  < ce^{-c_2(n-q)}N^{-\frac{\beta}{2}}.
\end{equation}

Let $\beta' = \frac{2p\beta}{\lambda}$ and assume $\beta$ small enough so that $\beta' < 1/4$.

For $q> 0$,  the walkers $(\xi_q^i)_{1\leq i\leq N}$ are independent conditionally to $\eta_{q-1}^N$ we have

$$[\eta_q^N - \phi_{q}(\eta_{q-1}^N)](f) = \eta_q^N\left (\frac{1}{N}\sum_{1\leq i \leq N} h_i\right ), $$
with $$h_i := f - S_{q-1, \eta_{q-1}^N}P_{q}(f)(\xi_{q-1}^i).
$$

Moreover, we have

$$\phi_{q}(\eta_{q-1}^N) = \frac{1}{N}\sum_{1\leq i\leq n} \mu_i ~~\text{with}~~ \mu_i =  \delta_{\xi_{q-1}^i}S_{q-1, \eta_{q-1}^N}P_{q}.$$

Since for any $i\in\llbracket 1, N \rrbracket$, we have $\mu_i(h_i) = 0$, we can apply Lemma 7.3.3 from \cite{bookFK}, and deduce that there exists $C\in\R$ such that

\begin{align*}
\E&\left [\left |\left [\eta_q^N - \phi_{q}(\eta_{q-1}^N)\right ]\left (F^N_{q, n}\right )\right |^{2\beta p}\right ]^{\frac{1}{2p}} \\
&= \E\left [\E\left (\left |\left [\eta_q^N - \phi_{q}(\eta_{q-1}^N)\right ]\left (F^N_{q, n}\right )\right |^{\lambda \beta'}~|~\xi_{q-1}\right )\right ]^{\frac{1}{2p}} \notag \\
&\leq \E\left [\E\left (\left |\left [\eta_q^N - \phi_{q}(\eta_{q-1}^N)\right ]\left (F^N_{q, n}\right )\right |^{\lambda}~|~\xi_{q-1} \right )^{\beta'}\right ]^{\frac{1}{2p}} \notag \\
&\leq \frac{C}{N^{\beta/2}}\E\left [\phi_{q}(\eta_{q-1}^N)\left (|F_{q, n}^N|^{\lambda}\right )^{\beta'} \right ]^{1/2p} .
\end{align*}

For $q=0$, the walkers are iid with common distribution $\eta_0$. The previous reasoning therefore holds using the convention $E(X ~|~ \xi_{-1}) = \E(X)$.

Applying the contraction property \eqref{contraction_h_process} with $\mu = \delta_x$ and $\eta = \delta_y$ we get the existence of $(a, \rho)\in\R^2_+$ such that

\begin{align}
|\Bar Q_{q, n}(f)(x) - &\Bar Q_{q, n}(f)(y)|  \notag \\
&\leq ae^{-\rho (n-q)}\left (1 + \frac{V(x)}{H_{q, n}(x)}\right )\left (1+\frac{V(y)}{H_{q, n}(y)}\right ). \label{ineq_diff_R}
\end{align}

By substituting \eqref{ineq_diff_R} into the definition of $F_{q, n}$ and applying Hölder's inequality along with Jensen's inequality, we obtain, for some $(a', \rho')\in\R_+^{*^2}$

\begin{align*}
&a'e^{\rho'  (n-q)}\E\left [\phi_{q}(\eta_{q-1}^N)\left (|F_{q, n}^N|^{\lambda}\right )^{\beta'} \right ]^{1/2p} \notag \\
\leq ~& 
\E\left [\phi_{q}(\eta_{q-1}^N)\{(H_{q, n}+V)^{\lambda}\}^{\beta'}\phi_{q}(\eta_{q-1}^N)\{(H_{q, n}+V)\}^{\lambda \beta'} \right. \\ 
&~~~~~~~~~~~~~~~~\times\left. \phi_{q}(\eta_{q-1}^N)(H_{q, n})^{-\lambda \beta'}\right ]^{\frac{1}{2p}} \notag \\
\leq ~&
\E\left [\phi_{q}(\eta_{q-1}^N)\{(H_{q, n}+V)^{\lambda}\}\right ]^{\frac{1}{4p}} \E\left [\phi_{q}(\eta_{q-1}^N)(H_{q, n})^{-2\lambda \beta'}\right ]^{\frac{1}{4p}}. 
\end{align*}

From our hypothesis on $Q_n$, we deduce from Lemma 3.2 in \cite{del2021stability}  that there exists a constant $c$ such that for any $(q', n')\in\N^2$, $H_{q', n'}\leq cV$. We can then conclude by choosing a small enough $\beta'$ and using Lemmas \ref{lmm:control_H_l} and  \ref{lmm:unif_bound_lyapu}. \qed \\

The proof of Theorem \ref{thm:optimal_predictor_unif_bound} is now relatively straightforward. \\

\textbf{Proof of Theorem \ref{thm:optimal_predictor_unif_bound}: } \\

Let $f\in\Ca_{V^\frac{\lambda}{4p}}(E)$. From the sub-additivity of the $L_p$-norm applied in \eqref{ont_step_decomp}, we have

\begin{align*}
E(|\eta_n^N(f)  -&\eta_n(f)|^p)^{\frac{1}{p}} \\ &\leq   \sum\limits_{q=0}^n \E\left ( \left | \phi_{q, n}(\eta_q^N)(f) - \phi_{q, n}(\phi_{q}(\eta_{q-1}^N))(f)\right |^p\right )^{\frac{1}{p}} .
\end{align*}

By applying Lemma \ref{lmm:bound_one_step}, we deduce that there exists $(c, \rho, \beta)\in\R\times\R^*_+\times(0,1]$ such that for any $N\in\N^*$ we have

\begin{equation*}
\E(|\eta_n^N(f)  -\eta_n(f)|^p)^{\frac{1}{p}} ~\leq~ \frac{c}{N^{\frac{\beta}{2}}}\sum\limits_{0\leq l \leq n} e^{-(n-l)\rho} 
~\leq~ \frac{c}{N^{\frac{\beta}{2}}(1-e^{-\rho})}.
\end{equation*}

This ends the proof of the theorem. \qed \\

\subsection{Ground state estimates}\label{subsection:grd_state_approx}

This subsection concentrates on proving Corollary \ref{cor:limit_bound}. We consider thus the time-homogeneous model.

Let $f\in\Ca_{V^\frac{\lambda}{4p}}(E)$. Notice that we can decompose the error made by the DMC method in the following way:

\begin{align*}
\E(|\eta_n^N(f) - &\eta_\infty(f)|^p)^\frac{1}{p} \\ &\leq \E(|\eta_n^N(f) - \eta_n(f)|^p)^\frac{1}{p} +|\eta_n(f) - \eta_\infty(f)|.
\end{align*}

Theorem \ref{thm:optimal_predictor_unif_bound} implies that there exists $(c_1, \beta)\in\R^*\times(0,1]$ such that

$$\sup_{n\in\N} \E(|\eta_n^N(f) - \eta_n(f)|^p)^\frac{1}{p} \leq \frac{c_1}{N^{\beta/2}}. $$

According to Theorem 4.3 in \cite{del2021stability}, there exists $(C_2, \omega)\in\R^{*^2}_+$ such that

$$ |\eta_n(f) - \eta_\infty(f)| \leq c_2e^{-\omega n}. $$

Hence

$$ \E(|\eta_n^N(f) - \eta_\infty(f)|^p)^\frac{1}{p} \leq \frac{c_1}{N^{\beta/2}} +  c_2e^{-\omega n}.$$

Thus, letting $a = \frac{1}{\omega}\ln(c_2/c_1)$ and $b = \frac{\beta}{2\omega}$, we deduce that there exists $c\in\R$ such that

\begin{equation*}
\sup_{n\ge a + b\ln(N)}\E(|\eta_n^N(f) - \eta_\infty(f)|^p)^\frac{1}{p} \leq \frac{c}{N^{\beta/2}}.
\end{equation*}

This concludes the proof. \qed

\section{Coupled harmonic oscillators}\label{section:harmonic_oscillator}

This section is devoted to a detailed study of the coupled harmonic oscillator. We recall that in this framework, we consider the flow of measures \eqref{limit-ga} starting from some initial measure $\eta_0$ on $\R^d$, and that the transition kernel and  the potential function considered are defined for some positive definite matrices $B$ and $S$ and some matrix $A$ by \eqref{def_P_A_B}.

This study will be carried out in four steps. First, we consider the general framework described at the beginning of Subsection \ref{subsection:harmonic_oscillator}. There, by performing a complete study of the Lyapunov functions associated with $P_{A, B}$ and $G_S$, we establish that under the simple condition $A^\prime SA < S$, the DMC method can be controlled uniformly in time.

Note that in the coupled harmonic oscillator framework, if the initial measure is normal, then the entire flow of measures consists of normal distributions. More specifically, the measures $\eta_n=\Na(m_n,\Omega_n)$ are defined with

$$
m_n=\Ea(\Omega_{n-1})~m_{n-1}\quad \mbox{\rm and}\quad \Omega_n=\Phi(\Omega_{n-1}),
$$
with the mappings
$$
\Ea(\Omega):=A(I+\Omega S)^{-1}\quad \mbox{\rm and}\quad\Phi(\Omega):=A (I+\Omega S)^{-1}\Omega A^{\prime }+B.
$$

Building on this, we will show that the condition $A^\prime SA < S$ is a stability condition ensuring that $m_n$ is a decreasing sequence converging to the null vector.

Returning to the standard harmonic oscillator, we then study a variation of the DMC method. This allows for uniform control of the method without any imposed condition on the matrices $A$ and $S$. This part also leads us to fully describe the ground state $h$ as well as the ground state energy associated with $P_{A, B}$ and $G_S$.

We continue with a result on the stability of the flow of measure $\eta_n$ with respect to the initial measure for the coupled harmonic oscillator. We prove that the harmonic oscillator tends to forget its initial condition exponentially fast. This result is of both theoretical and practical interest. It justifies that the constants in Corollary \ref{cor:conv_osci} should not be too dependent on the initial measure, and it proves that a bad initialisation of the walkers when running a simulation should have negligible consequences.
 
To complete our study, we examine the univariate harmonic oscillator. In this context, we prove that the expected error of the DMC method diverges to $\infty$ as soon as $A > 1$.  Additionally, we establish the asymptotic convergence  of this error for any value of $A$. Together, these results not only confirm the validity of our convergence theorem but also validate our choice to study the uniform-in-time convergence of the error rather than restricting our analysis to its asymptotic convergence.

\subsection{Lyapunov functions}\label{subsection:hamo_mutli_d}

This subsection is dedicated to the proof of Corollary \ref{cor:conv_osci}. Therefore, we place ourselves within the framework associated with this corollary.  We will conduct the proof assuming \eqref{Q(W) > W general}. If \eqref{Q(W) > W homogeneous} holds with $E_0<\Pbar$  then the demonstration is completely analogous.

Consider a family of Feller Markov transitions $(P_n)_{n\in\N}$  and a family of positive functions $(G_n)_{n\in\N}\in\Ca_0(\R^d)^\N$  bounded by $1$, verifying \eqref{upper_bound_density}, \eqref{lower_control} and \eqref{variation_G_GS}. It is clear
from Subsection \ref{subsection:markov_transition} that proving the existence of a continuous $P$-Lyapunov function $V \in \Ca_\infty(\R^d) $ makes \eqref{upper_bound_density} and \eqref{lower_control} sufficient conditions for \eqref{min-cond} to hold. To guarantee the existence of an appropriate $Q$-Lyapunov function, we need a result obtained by Kato in \cite{Kato:1966:PTL}. We present it here using the formulation provided in \cite{articleLiChiZhang}.

\begin{lmm}\label{lmm:continuity} Suppose that $D \subset \mathbb{R}$ is an interval, and let $A$ be a continuous function from $D$ to the space of real $d\times d$ matrices. In this case, there exist $d$ eigenvalues (counted with algebraic multiplicities)  of $A(t)$ which can be parameterized as continuous functions $\lambda_1(t)$, ..., $\lambda_d(t)$ from $D$ to $\R$. 
\end{lmm}

We can now ensure the existence of a $Q$-Lyapunov function under a simple matrix condition.

\begin{lmm}\label{lmm:hyp_1_gauss}
Assume that $A^TSA < S$. There exists $\alpha\in\R^*_+$ such that the function

\begin{equation}
V: x\in\R^d \mapsto \exp\left (\frac{\alpha}{2}x^TSx\right ),
\end{equation}

is a $Q$-Lyapunov function and it is integrable w.r.t $\eta_0$. 
\end{lmm}
\textbf{Proof : }\\

From \eqref{variation_G_GS}, we deduce that the r.h.s of \eqref{GQV} holds for $V$ and $G_n$. Then, together with
\eqref{upper_bound_density}, we deduce that it is enough to prove that $V$ is a Lyapunov function for $P_{A, B}$ with $\epsilon < 1/c$. Let's then compute $P_{A, B}(f)$ for any function of the form

$$f : x\in\R^d \mapsto \exp\left (\frac12 x^TFx\right ), $$

where $F$ is an invertible matrix such that $B^{-1} - F$ is positive definite.\\

In this setting, the Woodbury matrix identity provides the following equality:

\begin{equation*}
(B-F^{-1})^{-1}=B^{-1}-B^{-1}\left(B^{-1} - F\right)^{-1}B^{-1}.
\end{equation*}

We have then for any $x\in\R^d$, with $\Bar B_F := (B^{-1} - F)^{-1}$ :

\begin{align*}
&P_{A, B}(f)(x) \\
&= \frac {1}{(2\pi )^{d/2}\det(B)^{1/2}}\int_{\R^d} e^{-{\frac {1}{2}}\left [(Ax-y )^{\top }B ^{-1}(Ax-y ) - y^TFy\right ]}dy \notag\\
&\notag\\
&= f(x)~ \frac {\exp(\frac{1}{2}x^T(A^TB^{-1}\Bar B_F B^{-1} A-A^TB^{-1}A  - F)x)}{(2\pi )^{d/2}\det(B)^{1/2}} \notag \\
&~~~~~\times \int_{\R^d} \exp\left ({-{\frac {1}{2}}\left [(\Bar B_FB^{-1}Ax-y )^{\top }\Bar B_F ^{-1}(\Bar B_FB^{-1}Ax-y ) \right ]}\right )dy \notag
\end{align*}
This yields the formulae
\begin{equation}\label{M_f_calculations}
\begin{array}{l}
P_{A, B}(f)(x) \\
\\
\displaystyle= \sqrt{\frac {\det(\Bar B_F)}{\det(B)}}f(x)  \\ \\
~~~~~~~~~~~~\times \exp\left (\frac{1}{2}x^T(A^T[B^{-1}(B^{-1} - F)^{-1}B^{-1} -B^{-1}]A - F)x\right ) \\
\\
\displaystyle= \frac {1}{\sqrt{\det(I_d - BF)}}f(x) \exp\left (-\frac{1}{2}x^T(F +A^T( FB - I_d)^{-1}FA)x\right ).
\end{array}
\end{equation}

From those calculations, we deduce that $V$ is a Lyapunov function for $P_{A, B}$ if the matrices $B^{-1} - \alpha S$ and $S - A^T(I_d - \alpha SB)^{-1}SA$ are positive definite.

Let $\lambda_B$ be the greatest eigenvalue of $B$ and $\lambda_S$ be the greatest eigenvalue of $S$.  It is clear that for $\alpha \in (0 , \frac{1}{\lambda_B\lambda_S})$, $B^{-1} - \alpha S$ is positive definite.

Consider now the function

\begin{equation*}
\psi: \alpha \in \left [0 , \frac{1}{\lambda_B\lambda_S}\right ) \mapsto sp(S -A^T(I_d - \alpha SB)^{-1}SA)\in\R^d.
\end{equation*}

Here, $sp(M)$ represents the spectrum of a matrix $M$ with multiplicity taken into account.

Given the hypotheses on $A$ and $S$, we can conclude that $\psi(0) \subset \R_+^{*^d}$. Furthermore, by Lemma \ref{lmm:bound_one_step}, it is clear that $\psi$ is a continuous function. Since $\R_+^{*^d}$ is an open set, there exists $\bar{\alpha} \in \R$ such that for any $\alpha \in (0, \bar{\alpha})$, $\psi(\alpha) \subset \R_+^{*^d}$.

By choosing a sufficiently small value for $\alpha$ to ensure that $V$ is integrable w.r.t $\eta_0$, we can conclude. \qed \\

From this Lemma and the hypothesis on $G_n$, we deduce that the l.h.s of \eqref{GQV} holds as well. We can now focus on verifying that \eqref{Q(W) > W homogeneous} holds by proving the following Lemma.

\begin{lmm}\label{lmm:groud state oscilatore}
There exists a positive definite matrix $H$ such that
\begin{equation*}
\forall n\in\N, ~ Q_n(W) \ge \Pbar \times W\quad\mbox{with}\quad W(x) := \exp(-\frac{1}{2}x^THx)
\end{equation*}

\end{lmm}
\textbf{Proof : }\\

From \eqref{lower_control}, we have

\begin{equation*}
Q_n(W) \ge \frac{\Pbar} {E_{A', B', S}}\times G_S~P_{A', B'}(W).
\end{equation*}

Using \eqref{M_f_calculations}, we derive the following expression for $G_S(x)P_{A', B'}(W)(x)$ with $x\in\R^d$:

$$\frac {W(x)\exp\left (\frac{1}{2}x^T(H -{A'}^T( HB' + I_d)^{-1}HA' - S)x\right )}{\sqrt{\det(I_d + B'H)}} .
$$

Hence, chosing $H$ as the solution to the Riccati equation 
$$H-{A'}^T( HB' + I_d)^{-1}HA' - S = 0,$$
 we deduce that $$E_{A', B', S} = \frac {1}{\sqrt{\det(I_d + B'H)}}.$$ 
 
Thus: 
 
\begin{equation*}
Q_n(W) \ge \Pbar \times W.
\end{equation*} \qed \\

For $\alpha$ sufficiently small, $W^{-\alpha}$ is lower than $V$. The right-hand side of \eqref{Q(W) > W homogeneous} is then also verified.

Under the condition $A^TSA < S$, we have confirmed that all the assumptions concerning $P_n$, $G_n$, and $\eta_0$ in Theorem \ref{thm:optimal_predictor_unif_bound} are satisfied. We can therefore apply it to conclude on the proof of Corollary \ref{cor:conv_osci}.

Although justified by previous calculations, the meaning and implication of this matrix condition may not be immediately apparent. To remedy this, we consider the framework of the coupled harmonic oscillator to provide a property that illustrates how this condition is a stability condition that ensures that the mean of the flow of measures $\eta_n=\Na(m_n,\Omega_n)$ is a strictly decreasing sequence converging to zero.

\begin{prop}\label{prop:illustration_condition}
Let $B$ and $S$ be positive definite matrices  and $A$ be a matrix.
Consider a Gaussian distribution $\eta_0=\Na(m_0,\Omega_0)$ and consider the flow of measures $\eta_n=\Na(m_n,\Omega_n)$ given in (\ref{limit-ga}) with $P = P_{A, B}$ and $G = G_S$. The following implication holds

$$
A^\prime SA<S\quad \Longrightarrow \exists~C\in[0,1)~\text{s.t.}~ \forall ~k<n,~ \Vert m_n\Vert\leq C^{n-k}~~\Vert m_k\Vert
.$$
\end{prop}
\textbf{Proof : } \\

The flow of measures $\eta_n=\Na(m_n,\Omega_n)$ are defined with
$$
m_n=\Ea(\Omega_{n-1})~m_{n-1}\quad \mbox{\rm and}\quad \Omega_n=\Phi(\Omega_{n-1}),
$$
with the mappings
$$
\Ea(\Omega):=A(I+\Omega S)^{-1}\quad \mbox{\rm and}\quad\Phi(\Omega):=A (I+\Omega S)^{-1}\Omega A^{\prime }+B.
$$
We set
$$ \Ea_{k,n}(\Omega):=\Ea(\Phi^{n-1}(\Omega))~\ldots~\Ea(\Phi^{n-k}(\Omega))$$ 
as well as
$$\Ea_{n}(\Omega) = \Ea_{0,n}(\Omega)$$
and
$$\quad\Phi^n=\Phi\circ \Phi^{n-1}. $$

In this notation, we have
$$
m_n=\Ea_{k,n}(\Omega_0)~m_{k}.
$$
Note that
$$
A^\prime SA<S\quad \Longleftrightarrow\quad \overline{A}^\prime \overline{A}<I\quad\mbox{\rm with}\quad
\overline{A}:=S^{1/2}AS^{-1/2}.
$$
In the same vein, we have
$$
\Ea(\Omega)=S^{-1/2}~\overline{\Ea}(\Omega)~S^{1/2}$$
with
$$\overline{\Ea}(\Omega):=\overline{A}~(I+S^{1/2}~\Omega ~S^{1/2})^{-1}.
$$
Let $\Vert M\Vert_2:=\sqrt{\lambda_{max}(M^\prime M)}$ be the spectral norm, where $\lambda_{ max}(\Sa)$ stands for the maximal eigenvalue of a symmetric matrix $\Sa$. In this notation, we have the estimate
\begin{eqnarray*}
\Vert \overline{\Ea}(\Omega)\Vert_2&\leq &\Vert \overline{A}\Vert_2~\lambda_{max}((I+S^{1/2}~\Omega ~S^{1/2})^{-1}) \\
&\leq&\frac{\Vert \overline{A}\Vert_2}{\lambda_{\min}(I+S^{1/2}~\Omega ~S^{1/2})}\leq \Vert \overline{A}\Vert_2<1.
\end{eqnarray*}
Using the easily checked decomposition
$$
\Ea_{k,n}(\Omega)=S^{-1/2}\left(\overline{\Ea}(\Phi^{n-1}(\Omega))~\ldots~\overline{\Ea}(\Phi^{n-k}(\Omega))\right)S^{1/2},
$$
this yields the rather crude exponential decay
$$
\Vert \Ea_{k,n}(\Omega)\Vert_2\leq \Vert\overline{A}\Vert_2^{n -k}.
$$
We conclude that
$$
A^\prime SA<S\quad \Longrightarrow \Vert m_n\Vert\leq\Vert \overline{A}\Vert_2^{n-k}~~\Vert m_k\Vert\longrightarrow_{n\rightarrow\infty} 0
.$$

\subsection{Conditional free evolutions}\label{subsection:importance sampling}

This subsection focuses on the study of the importance sampling described in Section~\ref{subsection:DMC}, Section \ref{subsection:harmonic_oscillator} and on the proof of Corollary \ref{cor:optimal_filter}.

We consider the coupled harmonic oscillator, i.e, for some matrices $A$, $B$ and $S$, with $B$ and $S$ symmetric definite positive, we consider $P = P_{A, B}$ and $G = G_S$. Up to this point, we have established that the $L_p$-norm of the error made by the DMC method is uniformly bounded in time, with a convergence rate of $\frac{1}{N^{\beta/2}}$ for some $\beta \in(0, 1]$ when $A^TSA < A$.

Our aim is now to use Theorem \ref{thm:optimal_predictor_unif_bound} to prove that the approximation of the measures $\widehat \eta_n$ made by the DMC method - enhanced by the importance sampling scheme described in \eqref{approx_hat} - remains uniformly bounded in time, regardless of the value of $A$ and $S$. However, there's a trade-off involved: we will only have access to the measures at specific times. Indeed, despite the converging property that we are about to prove, the sequence of empirical measures $\widehat \eta_n^{(k)}$ only approximates the measures $\widehat \eta_l$ for $l\in k\N$.

Before proving the main result of this section, we conduction a short and necessaty study of the functions ${G}^{(k)}$.

\begin{lmm}\label{lmm:G chapeau n}
There exists a sequence of positive definite matrices $S_n$ such that

$$
{G}^{(n)}(x)\propto~
\exp{\left(-\frac{1}{2}~x^{\prime } S_nx\right)}.
$$
Moreover, there exists $ S_{\inf},  S_{\sup}$ symmetric positive definite such that $S_{\inf} <  S_n < S_{\sup}$.
\end{lmm}

\textbf{Proof : }\\
 
Using the same notations as the one used in the proof of Proposition \ref{prop:illustration_condition},
the Markov operators $\bar{Q}_{0,n} := {Q}_{0,n}/{Q}_{0,n}(1)$ are then given by

\begin{equation} \label{link A epsilon}
\delta_x\bar{Q}_{0, n}=\phi_{0,n}(\delta_x)=\Na\left(A_nx,  B_n\right) 
\end{equation}
with
$$A_n:= \Ea_n(0)
\quad \mbox{\rm and}\quad  B_n=\ \Phi^n(0).
$$

Using the product formula
$$
Q_{0, n}(1)(x)=\prod_{0\leq l<n}\phi_{0,l}(\delta_x)(G),
$$
we check that
$$
 Q_{0, n}(1)(x)\propto~
\exp{\left(-\frac{1}{2}~x^{\prime } S_nx\right)}
$$
with

\begin{align*}
    S_n:=&\sum_{0\leq l<n} \Ea_l(0)^{\prime }\left( S - S\left( \Phi^l(0)^{-1}+ S\right)^{-1} S\right) \Ea_l(0) \\
    =& \sum_{0\leq l<n} \Ea_l(0)^{\prime }\left( S^{-1} +  \Phi^l(0)\right)^{-1} \Ea_l(0)\geq 0 
\end{align*}

The matrices $ \Ea_l(0)^{\prime }\left( S^{-1} +  \Phi^l(0)\right)^{-1} \Ea_l(0)$ being positive for all $l\in\N^*$, $ S_n$ is lower bounded by $ S_1$.

Moreover, recalling that $S>0$ and using the results from Theorem~1.3. and Corollary~1.4 in ~\cite{del2022note} as well as the proof of 
Proposition 2.3, we know that there exists a positive definite matrix $P_\infty$, a matrix $E$ such that $\lambda_{max}(E) <1$ as well as a uniformly bounded sequence of matrices $L_n$  such that we have

$$  \Ea_n(0) = E^nL_n ~~\text{and}~~  \Phi^n(0) = P_\infty - E^n L_n P_\infty(E^n)^\prime $$

The sequence $ \Ea_l(0)^{\prime }\left( S- S\left( \Phi^l(0)^{-1}+ S\right)^{-1} S\right) \Ea_l(0)$ is then converging, thus there exists $C_1\in\R$ such that

$$  S_n \leq C_1\sum_{0\leq l<n} \Ea_l(0)^{\prime } \Ea_l(0)
$$

Since $E^n$ is converging to zero, there exists $n_1\in\N$  such that $\Vert E^{n_1}\Vert < 1$, with $\Vert . \Vert$   being the spectral norm.
Let $C_2 = C_1\times\max_{k\in\llbracket 0, n_1 \rrbracket} \Vert E^{k} \Vert \times \sup_{k\in\N} \Vert L_k \Vert$.
The spectral norm being sub-multiplicative, we have then

\begin{align*}
\Vert S_n \Vert \leq C_2\sum_{0\leq l<n}\left  \Vert E^{n_1\lfloor l/n_1 \rfloor}\right \Vert \leq n_1C_2\sum_{0\leq l<\lfloor n/n_1 \rfloor} \Vert E^{n_1} \Vert ^l.
\end{align*}

The norm of $S_n$ being bounded, we can conclude that there exists a positive definite matrix $ S_{\sup}$ such that 
$$\forall n\in\N, S_n \leq S_{\sup}$$.  \qed \\

The proof of Corollary \ref{cor:optimal_filter} is now relatively straightforward. \\

\textbf{Proof of Corollary \ref{cor:optimal_filter} : } \\
 
From Corollary \ref{cor:conv_osci}, Lemma \ref{lmm:G chapeau n} and \eqref{link A epsilon}, we deduce  that to prove the uniform convergence toward the Feynman-Kac measure of the DMC method associated with $ G^{(n)}$ and $ P^{(n)}$, it suffices to prove that

$$ A_n^\prime   S_n  A_n < S_n.$$

Let $\lambda_{ max}( S_n)$ be the greatest eigenvalue of $ S_n$. The matrix $ S - \lambda_{max}( S_n)I_d$ being negative, we get from Sylvester's law of inertia that for any non-singular matrix $ A$, $ A^\prime (S - \lambda_{max}( S_n)I_d) A$ is negative. By density of the invertible matrices, this holds for any matrix $ A$. Thus

$$
 A^{\prime }_n  S_n  A_n\leq \lambda_{max}( S_n) A_n^\prime  A_n.
$$

Denoting by $\Vert M\Vert_2$ be the spectral norm, we have then

$$
 A^{\prime }_n  S_n  A_n\leq \lambda_{max}( S_n)~\Vert  A_n\Vert_2^2~ I.
$$

From Lemma \ref{lmm:G chapeau n} and its proof, we deduce that $ A^{\prime }_n  S_n  A_n$ is converging to zero. Thus, for $n$ large enough, we have

$$ A_n^\prime   S_n  A_n < S_{\inf} <  S_n. $$

This concludes the proof. \qed \\

We conclude this section with an observation regarding the ground state energy associated to the coupled harmonic oscillator for any value of $A$, $B$ and $S$. Notice first that all of the hypotheses required for the existence of an eigen-triple \eqref{triplet-cond} are trivially met by taking $V : x\in\R^d \mapsto x^\prime x + 1$. Furthermore, with $x_0$ being the null vector, it is possible to express the triple as

$$
h(x) \propto h_\infty(x) := \exp\left (-\frac{1}{2}x^\prime S_{\infty}x\right ) ~~\text{,}~~ E_0 =  \frac{1}{\sqrt{\det(I_d + BS_\infty)}}$$ 
and
$$\eta_\infty\sim\mathcal{N}(0, P_\infty),
$$

with 
$$S_{\infty} := \lim\limits_{n \to \infty} S_n = \sum\limits_{l=0}^{\infty} \Ea_l(0)^{\prime }\left( S^{-1} +  \Phi^l(0)\right)^{-1} \Ea_l(0) $$ 
and 
$$P_\infty = \Phi^{\infty}(0).$$

Indeed, we have

\begin{align} \label{h_n}
h_n(x) &:= \exp(-\frac{1}{2}x^\prime S_{n}x) = \frac{Q_{0,n}(1)(x)}{Q_{0,n}(1)(x_0)} \notag \\
&\Rightarrow Q(h_n) =  \frac{Q_{0,n+1}(1)(x_0)}{Q_{0,n}(1)(x_0)} h_{n+1}.
\end{align}

Using the calcultions from  \eqref{M_f_calculations} we deduce that

\begin{align*}
Q_{0,n+1}(1)(x_0) &= Q[Q_{0,n}(1)](x_0) \\
&= Q_{0,n}(1)(x_0)\times Q(h_n)(x_0) \\
&= \frac{Q_{0,n}(1)(x_0)}{\sqrt{\det(I_d + BS_n)}}.
\end{align*}

Hence, going to the limit in \eqref{h_n} we obtain

$$ Q(h_\infty) = \frac{1}{\sqrt{\det(I_d + BS_\infty)}} h_\infty .$$

The value of $\eta_\infty$ directly follows from the computations of the flow of Gaussian measures $\eta_n$ in the proof of Proposition \ref{prop:illustration_condition}.

It is also worth noting that, by analogy with the calculations done in the proof of Lemma  \ref{lmm:groud state oscilatore} , $S_\infty$ is the solution to the Riccati equation

$$ S_\infty-{A}^\prime ( S_\infty B + I_d)^{-1}S_\infty A - S = 0. $$

\subsection{Stability w.r.t to the initial measure}

So far, in the context of the coupled harmonic oscillator, we have concentrated on examining the convergence of the DMC method with a fixed initial distribution. This naturally raises the question: does the choice of the initial measure significantly affect the control of the method ? To answer this question, we establish a stability result for the Feynman-Kac flow of measure \eqref{limit-ga} in the harmonic oscillator framework. Our analysis shows that the system tends to forget its initial condition exponentially fast. Consequently, regardless of the initial condition, the flow of the Feynman-Kac measure converges exponentially to the stationary measure $\eta_\infty$.

\begin{prop}
Let $S, B$ be two symmetric positive definite matrices on $\R^d$ and $A$ be a matrix. Consider $G = G_S$ and $P=P_{A,B}$. Let $\nu_1$ and $\nu_2$ be two probability measures on $\R^d$. There exists $(a, b)\in\R\times[0, 1)$ such that

$$ \Vert \phi_n(\nu_1) - \phi_n(\nu_2) \Vert_{tv}  \leq ab^n, $$

with the total variation distance between the measures $\nu_1$ and $\nu_2$ defined by
 
$$ \left\| \nu_1 - \nu_2 \right\|_{\text{tv}} := \sup_{\Vert f \Vert_{\infty} \leq 1} \left| \nu_1(f) - \nu_2(f) \right | 
$$
\end{prop}

\textbf{Proof : } \\

Let $f$ such that $\left\| f \right\|_{\infty} \leq 1$, then

\begin{align*}
&\phi_n(\nu_1)(f) - \phi_n(\nu_2)(f)  \\&= \frac{\nu_1 Q_n(f)}{\nu_1 Q_n(1)} - \frac{\nu_2 Q_n(f)}{\nu_2 Q_n(1)} \\
&= \psi_{Q_n(1)}(\nu_1)\Bar{Q}_n(f) - \psi_{Q_n(1)}(\nu_2)\Bar{Q}_n(f) \\
&= \int_{\R^d}\psi_{Q_n(1)}(\nu_1)(dx)\psi_{Q_n(1)}(\nu_2)(dy)\left [\Bar Q_n(f)(x) - \Bar Q_n(f)(y)\right ].
\end{align*}

We have then

\begin{align}\label{bound_integral_tv}
&| \phi_n(\nu_1)(f) - \phi_n(\nu_2)(f)| \notag \\ &\leq \int_{\R^d}\psi_{Q_n(1)}(\nu_1)(dx)\psi_{Q_n(1)}(\nu_2)(dy)\left \Vert \phi_n(\delta_x) - \phi_n(\delta_y)\right \Vert_{tv}.
\end{align}

From Pinsker's inequality, we have, with $Ent$ being the relative entropy - or Kullback–Leibler divergence

$$ \left \Vert \phi_n(\delta_x) - \phi_n(\delta_y)\right \Vert_{tv} \leq \sqrt{\frac{1}{2}~Ent(\phi_n(\delta_x)~|~\phi_n(\delta_y))}.$$

Using the results from \cite{Pardo:996837}, we know that for two Gaussian distributions $\mu_1 = \mathcal{N}(m_1, \Sigma_1^2)$ and $\mu_2 = \mathcal{N}(m_2, \Sigma_2^2)$, we have

\begin{align*}
Ent(\mu_1 | \mu_2) \notag = \frac{1}{2} [&log|\Sigma_2\Sigma_1^{-1}|  +  tr\left (\Sigma_2^{-1}\Sigma_1\right ) \\ & + \{(m_1 - m_2)^\prime \Sigma_2^{-1}(m_1 - m_2) - n].    
\end{align*}

Thus, using the notations from the proof of Property \ref{prop:illustration_condition}, $\phi_n(\delta_x)$ and $\phi_n(\delta_y)$ being Gaussians with the same covariance $\Omega_n$ and respective mean $\Ea_n(0)x$ and $\Ea_n(0)y$, we deduce that

$$Ent(\phi_n(\delta_x)~|~\phi_n(\delta_y)) = (x - y)^\prime \Ea_n(0)^\prime ~\Omega_n~\Ea_n(0)(x-y).$$

From Theorem 1.3. and Corollary 1.4 in ~\cite{del2022note}, we know that $\Omega_n$ is uniformly bounded in time and that $\Ea_n(0)$ converges exponentially fast to zero. Thus, there exists $(a, b)\in\R\times[0, 1)$ such that

$$Ent(\phi_n(\delta_x)~|~\phi_n(\delta_y)) \leq ab^n\Vert x - y\Vert^2.$$

Combining this with Pinsker's inequality and \eqref{bound_integral_tv}, we get, for some $(a, b)\in\R\times[0, 1)$

\begin{align*}
| \phi_n(\nu_1)(f) - \phi_n(\nu_2)(f)| \leq ab^n&[\psi_{Q_n(1)}(\nu_1)(||\cdot||) \\ &+ \psi_{Q_n(1)}(\nu_2)(||\cdot||)],
\end{align*}

where $||\cdot||$ represents the function that maps $x \in \mathbb{R}^d$ to its norm.

Using the computations from Lemma \ref{lmm:G chapeau n}, we know that

$$Q_n(1) \propto~
\exp{\left(-\frac{1}{2}~x^{\prime } S_nx\right)} := q_n, $$

and that there exists two symmetric positive definite matrices $S_{\sup}$ and $S_{\inf}$ such that for any $n\in\N$, $S_{\inf}\leq S_n \leq S_{\sup}$.

Considering the function $q_{\inf} := \exp{\left(-\frac{1}{2}~x^{\prime } S_{\inf}x\right)}$ and  $q_{\sup} := \exp{\left(-\frac{1}{2}~x^{\prime } S_{\sup}x\right)}$, we have then

$$\Psi_{Q_n(1)}(\nu_1)(||\cdot||) =  \frac{\nu_1(||\cdot||q_n)}{\nu_1(q_n)} \leq \frac{\nu_1(||\cdot||q_{\inf})}{\nu_1(q_{\sup})}. $$

The function $||\cdot||q_{\inf}$ and $q_{\sup}$ being bounded, we deduce that the previous upper bound is finite.

Hence, for any $f$ such that $\left\| f \right\|_{\infty} \leq 1$, we have

$$ | \phi_n(\nu_1)(f) - \phi_n(\nu_2)(f)| \leq ab^n\left [\frac{\nu_1(||\cdot||q_{\inf})}{\nu_1(q_{\sup})}+  \frac{\nu_2(||\cdot||q_{\inf})}{\nu_2(q_{\sup})} \right ].$$

We can thus conclude. \qed \\

\subsection{Divergence  and fluctuation estimates}\label{subsection:hamo_1_d}

In the previous subsections, we presented a simple sufficient condition for controlling the DMC method and introduced an importance sampling technique that satisfies this criterion. However, it is natural to question the robustness of this condition and whether it is necessary to use importance sampling. Specifically, for the uni-dimensional harmonic oscillator, the convergence condition reduced to $A^2 < 1$, and we will prove divergence of the DMC method when this stability condition is not met.

Within this framework, we can break down the evolution of the walkers into two distinct steps, a mutation transition and a selection transition

\begin{align*}
\left (\xi_n^{i}\right )_{i\in\llbracket 1, N \rrbracket} \in \R^N &~~\xrightarrow{selection}~~ \left (\widehat \xi_n^{i}\right )_{i\in\llbracket 1, N \rrbracket} \in \R^N \\
&~~~\xrightarrow{mutation} \left (\xi^{i}_{n+1}\right )_{i\in\llbracket 1, N \rrbracket}  .
\end{align*}
 
The initial configuration $\left(\xi_0^{i}\right)_{i\in\llbracket 1, N \rrbracket}$ is determined by sampling $N$ independent random variables from the distribution $\eta_0$. The selection transition involves the sampling of $N$ independent random variables $\left(\widehat \xi_n^{i}\right)_{i\in\llbracket 1, N \rrbracket}$ using the weighted distributions

\begin{align*}
\epsilon_n(\eta_n^N)&G_S(\xi^i_n)\delta_{\xi_n^{i}} \\ &+ (1 - \epsilon_n(\eta_n^N)G_S(\xi^i_n))\sum_{k\in\llbracket 1, N \rrbracket}\frac{e^{-\frac{S}{2}\xi_n^{k^2}}}{\sum\limits_{j\in\llbracket 1, N \rrbracket} e^{-\frac{S}{2}\xi_n^{j^2}}}\delta_{\xi_n^{k}}.
\end{align*}

The mutation transition is defined using a family of Gaussian random variables with zero-mean and unit variance $(V^i_n)_{i\in\llbracket 1, N \rrbracket}$ such that

\begin{equation*}
\xi_{n}^{i} = A\widehat \xi_{n-1}^{i} + \sqrt{B}V^i_n.
\end{equation*}

The measures $(\eta_n)_{n\in\N}$ can be described exhaustively  using the Kalman filter equations. It provides us with the mean and variances $(m_n, \sigma_n^2)$ of the Gaussian random variables $\eta_n$ with the recurrent equations

\begin{equation}
\left\{
    \begin{array}{ll}
        \displaystyle   m_{n+1} =  \frac A{1+S\sigma_n^2}m_n \\
        \\
        \displaystyle \sigma_{n+1}^2 = \frac{A^2\sigma_n^2}{1+S\sigma_n^2} + B.
    \end{array}
\right. \label{Kalman}
\end{equation}

In this scenario, when the condition $A^2 > 1$ is met, it is possible to prove that the DMC's error in approximating the Feynman-Kac measure does not admit a uniform-in-time bound. It is properly stated in Property \ref{thm:div_oscillo}, and we can now conduct its proof.\\

\textbf{Proof of Proposition \ref{thm:div_oscillo}}: \\

For any $n\ge2$, we know from \eqref{Kalman} that

\begin{equation*}
\eta_n(I) := m_n = \frac {A^2}{(1+S\sigma_{n-1}^2)(1+S\sigma_{n-2}^2)}m_{n-2} \leq A^2 m_{n-2} .
\end{equation*}

For any $n\in\N^*$, let $\xi^*_n = \min\limits_{i\in\llbracket 1, N \rrbracket} \xi^i_n$ and define the random variables $V^*_n$ in the following way:

\begin{equation*}
V^*_n =  
        -\max\limits_{i\in\llbracket 1, N \rrbracket} \left |V^i_n \right | .
\end{equation*}

By definition of the evolution of the walkers, there exits $(i,j)\in\llbracket 1, N \rrbracket $ such that

\begin{align*} 
\xi^*_{2n} &= A^2\xi^j_{2n-2} + \sqrt{B}V_{2n-1}^i + A\sqrt{B}V_{2n-2}^j \\
&\ge  A^2\xi^*_{2n-2} + \sqrt{B}V_{2n-1}^* + |A|\sqrt{B}V_{2n-2}^*.
\end{align*}

Thus

\begin{align*}
\eta_{2n}^N(I) - \eta_{2n}(I) \ge A^2(\xi^*_{2(n-1)}& - m_{2(n-1)}) \\ + &\sqrt{B}V_{2n -1}^*  + |A|\sqrt{B}V_{2n-2}^* .
\end{align*}

Iterating the process, we obtain

\begin{align*}
\frac{\eta_{2n}^N(I) - \eta_{2n}(I)}{A^{2n}} \ge (\xi_0^* - &m_0) + \sqrt{B}\sum_{1\leq k\leq n}\frac{V_{2k-1}^*}{A^{2k}} \\ &+ |A|\sqrt{B}\sum_{1\leq k\leq n}\frac{V_{2(k-1)}^*}{A^{2k}} .
\end{align*}

For any sequence of $N$ independent centred Gaussian random variables $U_i$ with unit variance, we have

$$ \E\left [\max\limits_{1\leq i\leq N} |U_i|\right ] \leq \sqrt{2\log(2N)}. $$

This inequality is obtained by using Jensen's inequality as follows, with $t = \sqrt{2\log(2N)}$

\begin{align*}
    \exp\left [t\E\left (\max\limits_{1\leq i\leq N} |U_i|\right )\right ]&\leq \E\left [\exp\left (t\max\limits_{1\leq i\leq N} |U_i|\right )\right ] \\ &\leq \sum_{i=1}^N\E\left [\exp\left (t|U_i|\right )\right ],
\end{align*}

and noticing that

\begin{align*}
\E\left [\exp\left (t|U_i|\right )\right ] &= 2\int_0^{+\infty}\exp\left (-\frac{(x - t)^2 + t^2}{2}\right )dx \\ &\leq 2 \exp(t^2/2).
\end{align*}

Then, on the event

\begin{equation}
\Omega_\epsilon := \left \{ \xi_0^* \ge \epsilon + m_0 +  \frac{2\sqrt{B}(1 + |A|)}{A^2-1}\sqrt{2\log(2N)}\right \},
\end{equation}

with $\epsilon\in\R^*_+$, we have

\begin{equation}
\E[\eta_{2n}^N(I) - \eta_{2n}(I) | \xi_0^*]  ~\ge~ \epsilon A^{2n} ~\xrightarrow{n\rightarrow +\infty}~ +\infty.
\end{equation}

Integrating over $\xi_0^*$ we deduce

\begin{equation}
\E[|\eta_{2n}^N(I) -\eta_{2n}(I)| |] ~\ge~ \epsilon A^{2n}\Prb(\Omega_\epsilon).
\end{equation}

We can then conclude by noticing

\begin{align*}
\Prb(\Omega_\epsilon) &= \eta_0\left \{ \left [\epsilon + X_0 + \frac{2\sqrt{B}(1+|A|)}{A^2-1}\sqrt{2\log(2N)} , +\infty \right )\right \}^N \\ &> 0.
\end{align*}  \qed \\

For the case \( A = 1 \), we are unable to determine whether a uniform bound exists. To the best of our knowledge, the strongest divergence-type result established so far is a linear bound on the variance of the unnormalized measure for \( R = S = 1 \) obtained in \cite{WHITELEY20121840}. In that work, the authors employed a stronger exponential drift condition compared to our linear condition \eqref{lyap-cond}. As a result, they used a quadratic Lyapunov function, whereas we adopted an exponential Lyapunov function, which enabled us to apply Corollary~\ref{cor:conv_osci} to a broader range of test functions.

The divergence result in Proposition \ref{thm:div_oscillo} highlights the importance of studying non-asymptotic uniform convergence results rather than relying solely on central limit theorems (CLTs). Extensive research has been devoted to CLTs and, under appropriate assumptions. We quote the first studies in this field \cite{del1999central,DelMoral2000} mainly based on uniformly bounded potential and test functions.

More general fluctuation theorems that apply to more general models including diffusion-type processes with Lipschitz drift and diffusion functions as well as test functions with at most polynomial growth  are discussed in
 \cite{delmojacod,delmojacod2}. We can formulate the following result :

\begin{lmm}\label{lmm:TCL }
In the context of proportional selection/reconfiguration, for any $n\ge 1$, we have the following convergence in law as N tends toward $+\infty$

\begin{equation*}
\sqrt{N}[\eta_n(I) - \eta_n^N(I)] \xrightarrow[N \rightarrow \infty]{\mathcal{L}} \mathcal{N}(0 , \sigma_n^2) ,
\end{equation*}
with the asymptotic variance
\begin{equation*}\sigma_n^2 := \sum\limits_{p=0}^n \eta_p\left [\left (\frac{Q_{p,n}(1)}{\eta_pQ_{p,n}(1)} \Bar{Q}_{p,n}(I - \eta_n(I))\right )^2\right ] .
\end{equation*}
\end{lmm}

Related asymptotic variance formulae and comparisons are discussed in~\cite{moral2018note} in the context of random walks with absorbing barriers, including geometric killing rates and local reflection moves.

In \cite{10.1214/12-AAP909}, Nick Witheley also obtains a uniform time bound on the asymptotic variance. However, our next proposition shows that, for a given system, the DMC method can have a uniform time bound on its asymptotic variance, despite the fact that its non-asymptotic variance is unbounded.

\begin{prop}\label{prop-fi}
Let $G : x\in\R \mapsto e^{-\frac{x^2}{2}S}$, and $P$ such that $\delta_x P \sim \mathcal{N}(Ax, B)$ for some $(A, B, S) \in (1, +\infty]\times\R_+^{*^2}$. 
There exists an intial distibution $\eta_0$ such that

\begin{equation}
\sup_{n\in\N} \sigma^2_n < \infty ~~\text{and}~~ \sup_{n\in\N} \E\left (|\eta_n(I) - \eta_n^N(I)|^p\right )^{\frac{1}{p}} = +\infty, \label{div_conv}
\end{equation}

where $\sigma^2_n$ is defined as in Lemma \ref{lmm:TCL }.
\end{prop}

\textbf{Proof : } \\

Considering $\eta_0 \sim\mathcal{N}(m_0, \sigma_0^2)$ and using the Kalman filter's equations, we are able to fully describe the measures $\eta_n$. For $n\in\N^*$ we have

\begin{equation*}
\eta_n \sim\mathcal{N}(m_n, \sigma_n^2) ~~\text{with}~~ 
 \left\{
    \begin{array}{ll}
        \displaystyle m_{n+1} = \frac{A}{1 + S\sigma_n^2}~m_n \\
        \\
      \displaystyle  	\sigma_{n+1}^2 = \frac{A^2}{1+S\sigma_n^2} + B.
    \end{array}
\right.
\end{equation*}

The limit measure is given by $\eta_\infty = \mathcal{N}(0, \sigma_\infty^2)$ where $\sigma_\infty^2$ is the fixed point of the function $x\in\R_+\mapsto \frac{A^2}{1+Sx} + B$. From now on, we assume that $\eta_0 = \eta_\infty$. We have then

\begin{align*}
\sigma_n^2 &= \sum_{p=0}^n  \frac{\eta_\infty\left (Q_{p,n}(I)^2\right )}{\left (\eta_\infty Q_{p,n}(1)\right )^2}.
\end{align*}

Using the calculations and notations from Lemma \ref{Q_1_Q_I}, we have then

\begin{align*}
    \displaystyle    \eta_\infty\left (Q_{p,n}(1)\right ) &=   \displaystyle\frac {\lambda_{n-p}}{\sqrt{2\pi \sigma_\infty^2}}\int_{\R} e^{-{\frac {y^2}{2\sigma_\infty^2}}} e^{-\frac{y^2}{2}q_{n-p}}dy  \\ &= \frac{\lambda_{n-p}}{\sqrt{q_{n-p}\sigma_\infty^2 + 1}},
\end{align*}

as well as

\begin{align*}
    \displaystyle  \eta_\infty\left (Q_{p,n}(I)^2\right ) &=   \displaystyle \frac{\mu_{n-p}^2}{\sqrt{2\pi \sigma_\infty^2}}\int_{\R} y^2e^{-{\frac {y^2}{2\sigma_\infty^2}}} e^{-y^2q_{n-p}}dy  \\ &= \frac{\mu_{n-p}^2\sigma_\infty^2}{(2q_{n-p}\sigma_\infty^2 + 1)^{3/2}}.
\end{align*}

Then

\begin{align*}
\sigma_n^2 &=  \sigma_\infty^2 \sum_{p=0}^n  \frac{\mu_{n-p}^2}{\lambda_{n-p}^2} \frac{(q_{n-p}\sigma_\infty^2 + 1)^{3/2}}{(2q_{n-p}\sigma_\infty^2 + 1)^{3/2}}.
\end{align*}

As a solution to a Riccati equation, $q_n$ converges towards some $q_\infty\in\R^*_+$ such that

\begin{equation*}
q_\infty = \frac{A^2q_\infty}{1+q_\infty B} + S.
\end{equation*}

Since $A>1$, we obtain that $1 + q_\infty B =  A^2 + S/q_\infty> A$. Thus, there exists $C\in (0,1)$ such that, for $n$ large enough, we have

$$ \mu_{n+1} \leq C\frac{\mu_n}{\sqrt{1 + q_nB}}. $$

Thus, comparing the definition of $\lambda_n$ and the previous bound, we deduce that there exists $\alpha \in \R$ such that we have

\begin{equation*}
\frac{\mu_{n-p}^2}{\lambda_{n-p}^2} \frac{(q_{n-p}\sigma_\infty^2 + 1)^{3/2}}{(2q_{n-p}\sigma_\infty^2 + 1)^{3/2}} \leq \alpha C^{2(n-p)}.
\end{equation*}

We can thus conclude on the right-hand side. of \eqref{div_conv}. The left-hand side is a direct application of Proposition \ref{thm:div_oscillo}. \qed

\section*{Acknowledgments} MC thanks the European Research Council (ERC)
under the European Union’s Horizon 2020 research and
innovation program (grant agreement no. 863481) for
financial support.

Luc de Montella would like to thank Naval Group for their financial support through the CIFRE PhD program.

We also thank the anonymous reviewer for her/his excellent suggestions for improving the paper. Her/his detailed comments greatly improved the
presentation of the article.
\section*{Appendix}

\subsection*{Some technical Lemmas}

\begin{lmm}\label{lmm:remove_gibbs}
Let $(G, V)\in\Ca(E)\times\Ca_\infty(E)$ be positive functions and  $K\subset E$ be such that the right-hand side of \eqref{GQV} holds. There exists $c\in\R_+$ such that for any probability measure $\mu$ on $E$ :

\begin{equation*}
\psi_G(\mu)(V) \leq \mu(V) + c \label{ineq_psi}.
\end{equation*} 
\end{lmm}

\textbf{Proof :}
Let $\Bar V := \mathbb{1}_{\R\setminus K} V$. For any $(x,t)\in E^2$, we have:

\begin{align*}
0\ &ge (G(x) - G(y))(\Bar V(x) - \Bar V(y)) \\ &= G(x) \Bar V(x) + G(y) \Bar V(y) - G(y) \Bar V(x) - G(x) \Bar V(y).
\end{align*}

Integrating with respect to the probability measure $\mu$ over both $x$ and $y$, we obtain:

\begin{equation*}
\mu(G \Bar V)  - \mu(G)\mu( \Bar V) \leq 0.
\end{equation*}

Hence

\begin{equation*}
\psi_G(\mu)( \Bar V) - \mu( \Bar V) = \frac{\mu(G \Bar V) - \mu(G)\mu( \Bar V)}{\mu(G)} \leq 0.
\end{equation*}

Thus

\begin{align*}
\psi_G(\mu)(V) \leq \psi_G(\mu)(\bar V) + \sup_K V &\leq \mu(\bar V) + \sup_K V \\ &\leq  \mu(V) + \sup_K V.
\end{align*} \qed \\

\begin{lmm}\label{lmm:unif_bound_lyapu}

Let $V\in\Ca_\infty(E)$ be a $Q$-Lyapunov function, then

\begin{equation}
\left\{
    \begin{array}{ll}
        \sup_{\substack{(q,n, N)\in\N^3 \\ q \leq n}}\E\left [ \phi_{q , n}(\eta_q^N)(V)\right ] < +\infty \\ \\
        \sup_{\substack{(q,n)\in\N^2 \\ q \leq n}}  \phi_{q, n}(\eta_q)(V) < +\infty.
    \end{array}
\right.  
\label{unif_bound_phi_l}
\end{equation}  
\end{lmm}

\textbf{Proof: } 
First, let's prove that for any $\gamma \in \Pa_V(E)$, we have, with $(\epsilon, c)\in [0,1)\times\R$ defined in~\eqref{lyap-cond},

\begin{equation}
\sup_{\substack{(q,n)\in\N^2 \\ q \leq n}}~~\phi_{q, n}(\gamma)(V) \leq \gamma(V) + \frac{c}{1 - \epsilon}. \label{remove_phi_2}
\end{equation}

To do so, we begin by using the $Q$-Lyapunov property of $V$ as well as Lemma \ref{lmm:remove_gibbs} for some $l\in\N^*$ to deduce that

\begin{align*}
\phi_{q, n}(\gamma)(V)  &= \frac{\phi_{q, n-1}(\gamma)(G_nP_n(V))}{\phi_{q, n-1}(\gamma)(G_n)} 
\leq \epsilon~ \phi_{q, n-1}(\gamma)(V) + c'.
\end{align*}

By iterating the process, we obtain \eqref{remove_phi_2}. We now prove that

\begin{equation}
\sup_{n\in\N}~~\E\left [ \eta_n^N(V)\right ] \leq \eta_0(V) + \frac{c'}{1 - \epsilon}. \label{remove_phi_3}
\end{equation}

Notice first that

\begin{align*}
\E\left [ \eta_n^N(V) ~|~\eta_{n-1}^N\right ] &= \frac{1}{N}\sum_{i=1}^N\eta_{n-1}^NS_{n-1, \eta_{n-1}^N}P_n(V) \\ &= \psi_{G_{n-1}}(\eta_{n-1}^N)(P_n(V)) .
\end{align*}

Then, using as previously the $Q$-Lyapunov property of $V$ and Lemma \ref{lmm:remove_gibbs} we obtain

\begin{align*}
\psi_{G_{n-1}}(\eta_{n-1}^N)(P_n(V)) 
\leq \epsilon~ \eta_{n-1}^N(V) + c'.
\end{align*}

By iterating the process, we obtain

\begin{equation}
\E\left [ \eta_n^N(V)\right ] \leq \epsilon^n \eta_0(V) + c'\sum_{i=0}^{n-1}\epsilon^i \notag \leq \eta_0(V) + \frac{c'}{1-\epsilon}. \label{remove_empirical}
\end{equation}

Thus, by combining \eqref{remove_phi_2} and \eqref{remove_phi_3}, we can conclude regarding the first part of \eqref{unif_bound_phi_l}. The second part is obtained by proceeding in a strictly analogous way. \qed \\

\subsection*{Proof of Lemma \ref{lmm:control_H_l} }
We first consider the case where only \eqref{Q(W) > W general} holds and prove that

\begin{equation}
\phi_{0, q}(\gamma)( Q_{q, n}(1)) \leq C\Pbar^{n-q-1}. \label{boung_Q_n_k}
\end{equation}

For $q\leq n+2$, we have

\begin{align*}
Q_{q, n}(1)(x) &= Q_{q, n-2}(G_{n-1}P_{n}(G_n))(x) \\
&\leq \Pbar ~Q_{q, n-2}(G_{n-1})(x) \leq \Pbar~ Q_{q, n-1}(1)(x).
\end{align*}

Iterating the process, we deduce

$$Q_{q, n}(1) \leq \Pbar^{n-q-1} G_{q-1}. $$

Using \eqref{control_phi_G} we deduce \eqref{boung_Q_n_k}.

For $q\ge 1$ and $\mu\in\{\eta_q^N, \phi_{q}(\eta_{q-1}^N)\}$, we have then

\begin{align*}
\frac{1}{\mu(H_{q, n}^{\gamma})} &= \phi_{q, n}(\mu)(W)\times\frac{\phi_{0, q}(\gamma) Q_{q, n}(1)}{\mu(Q_{q, n}(W))} \\ &\leq C\frac{\Pbar^{n-q-1} \phi_{q, n}(\mu)(W)}{\mu(Q_{q, n}(W))}. 
\end{align*}

From Holder's inequality, Jensen's inequality and the hypothesis on $W$, we get :

\begin{align*}
\E\left [\mu(H_{q, n})^{-\beta}\right ] \notag 
\leq ~~&C^\beta\E\left [\phi_{q, n}(\mu)(W)^{2\beta}\right ]^{\frac{1}{2}}\\ &\times\E\left [\Pbar^{2\beta(n-q-1)}\mu(Q_{q, n}(W))^{-2\beta}\right ]^{\frac{1}{2}} \notag \\
\leq ~&\frac{C^{\beta}}{\Pbar^\beta}\E\left [\phi_{q, n}(\mu)(V)^{2\beta}\right ]^{\frac{1}{2}}\E\left [\mu(W^{-2\beta})\right ]^{\frac{1}{2}} .
\end{align*}

Then, by choosing $\beta$ small enough such that $W^{-2\beta} \leq \Bar W$ and $V^{2\beta} \leq \Bar V$, where $\Bar W$  and $\Bar V$ are $Q$-Lyapunov functions, we obtain

\begin{align*}
\sup_{\substack{(q,n, N)\in\N^3 \\ q \leq n}} \E\left [\mu(H_{q, n}^\gamma)^{-\beta}\right ] 
\leq \frac{C^{\beta}}{\Pbar^\beta}&\sup_{\substack{(q,n, N)\in\N^3 \\ q \leq n}} \E\left [\mu(\Bar V)\right ]^{\frac{1}{2}} \\ &\times \sup_{\substack{(q,n, N)\in\N^3 \\ q \leq n}} \E\left [\mu(\Bar W)\right ]^{\frac{1}{2}}.
\end{align*}

We can conclude using Lemma \ref{lmm:unif_bound_lyapu}.

The demostration for time-homogeneous models with $E_0 < \Pbar$ is analogous. Indeed the equivalent of \eqref{boung_Q_n_k} is obtained by noticing that

\begin{align*}
    \phi_{0, q}(\eta_\infty) Q_{q, n}(1) = \eta_\infty Q_{q, n-1}(G) &= \eta_\infty Q_{q, n-1}(1) \eta_\infty(G) \\ &= E_0 \times\eta_\infty Q_{q, n-1}(1) = E_0^{n-q}.
\end{align*}
The rest of the proof follows the same arguments, thus it is skipped. This ends the proof of the lemma.\qed \\

\begin{lmm}\label{Q_1_Q_I}
Let $G : x\in\R \mapsto e^{-\frac{x^2}{2}S}$, and $P$ such that $\delta_x P \sim \mathcal{N}(Ax, B)$ for some $(A, B, S) \in \R\times\R_+^{*^2}$. Denoting by $Q^n$ the $n$-times composition of the operator $Q$, we have

\begin{equation}
Q^n(1)(x) = \lambda_n e^{-\frac{x^2}{2}q_n} ~~\text{and}~~ Q^n(I)(x)= \mu_n~x~e^{-\frac{x^2}{2}q_n},
\end{equation}
with the parameters $(q_0,\lambda_0,\mu_0)=(S,1,A)$ and
\begin{align*}
 \left\{
    \begin{array}{ll}
        q_{n+1} = \displaystyle \frac{A^2q_n}{1 + q_nB} +  S\\
        \\
        \lambda_{n+1} = \displaystyle \frac{\lambda_n}{\sqrt{1 + q_nB}} \quad \mbox{and}\quad
        \mu_{n+1} = \displaystyle \frac{A\mu_n}{(1 + q_nB)^{3/2}}.
    \end{array}
\right.
\end{align*}
\end{lmm}

\textbf{Proof } :

For $n=0$ the result is immediate. Assume that it holds for some $n\in\N$. In this situation, we have

\begin{align*}
Q^{n+1}(1)(x) &= Q(Q^n(1))(x) = \lambda_n Q\left [y\mapsto e^{-\frac{y^2}{2}q_n}\right ](x) \\
&= \frac{\lambda_n}{\sqrt{2\pi B}} e^{-\frac{x^2}{2}S}\int_\R e^{-\frac{y^2}{2}q_n -\frac{(Ax - y)^2}{2B}}dy \\
	&= \frac{\lambda_n}{\sqrt{2\pi B}} e^{-\frac{x^2}{2}\left (\frac{A^2}{B} - \frac{A^2}{B(1 + q_nB)} +  S\right )}\\
 &~~~~~~~~~~~~~~~~~~\times\int_\R e^{-\frac{(1 + q_nB)\left (\frac{A}{1 + q_nB}x - y\right )^2}{2B}}dy \\
	&= \frac{\lambda_n}{\sqrt{1 + q_nB}} e^{-\frac{x^2}{2}\left (\frac{A^2q_n}{1 + q_nB} +  S\right )}.
\end{align*}

Similarly, we have

\begin{align*}
Q^{n+1}(I)(x) &= Q(Q^n(I))(x) = \mu_n Q\left [y\mapsto ye^{-\frac{y^2}{2}q_n}\right ](x) \\
	&= \frac{\mu_n}{\sqrt{2\pi B}} e^{-\frac{x^2}{2}\left (\frac{A^2q_n}{1 + q_nB} +  S\right )} \\ &~~~~~~~~~~~~~~~~~~\times\int_\R ye^{-\frac{(1 + q_nB)\left (\frac{A}{1 + q_nB}x - y\right )^2}{2B}}dy \\
	&= \frac{A\mu_n}{(1 + q_nB)^{3/2}} x e^{-\frac{x^2}{2}\left (\frac{A^2q_n}{1 + q_nB} +  S\right )}.
\end{align*}

This ends the proof. \qed

\bibliographystyle{unsrt}
\bibliography{main}{}
\bibliographystyle{plain}

\end{document}